\documentclass[11pt]{amsart}

\usepackage{fullpage, graphicx, enumitem}
\usepackage{url}
\usepackage{verbatim}
\usepackage{bbm}
\usepackage{amssymb,amsmath}
\usepackage{amsfonts,savesym,graphicx,bm,amsthm}
\usepackage{color}
\usepackage[mathscr]{eucal}
\usepackage[all]{xy}
\usepackage{tikz}
\usepackage{placeins}

\definecolor{hot}{RGB}{65,105,225}
\usepackage[pagebackref=true,colorlinks=true, linkcolor=hot ,  citecolor=hot, urlcolor=hot]{hyperref}

\newtheorem{theorem}{Theorem}[section]
\newtheorem{lemma}[theorem]{Lemma}
\newtheorem{conjecture}[theorem]{Conjecture}
\newtheorem{theorem-definition}[theorem]{Theorem-Definition}

\newtheorem{proposition}[theorem]{Proposition}

\newtheorem{definition-theorem}[theorem]{Definition-Theorem}

\newtheorem{corollary-definition}[theorem]{Corollary-Definition}
\newtheorem{definition-proposition}[theorem]{Definition-Proposition}

\theoremstyle{definition}

\newtheorem{definition}[theorem]{Definition}

\newtheorem{remark}[theorem]{Remark}

\newtheorem{subs}[theorem]{ }

\numberwithin{equation}{section}

\def\bC{\mathbb{C}}
\def\be{\begin{equation}}
	\def\ee{\end{equation}}

\author{Yifan Chen}
\address{Department of Mathematical Sciences, Tsinghua University, Beijing, 100084, P. R. China.}
\email{c-yf20@tsinghua.org.cn}

\author{Quan Shi}
\address{Department of Mathematical Sciences, Tsinghua University, Beijing, 100084, P. R. China.}
\email{shiq24@mails.tsinghua.edu.cn / thusq20@gmail.com}

\author{Yongxin Xu}
\address{Zhili College, Tsinghua University, Beijing, 100084, P. R. China.}
\email{xuyongxi@mails.tsinghua.edu.cn}

\author{Huaiqing Zuo}
\address{Department of Mathematical Sciences, Tsinghua University, Beijing, 100084, P. R. China.}
\email{hqzuo@mail.tsinghua.edu.cn}

\title{On the monodromy conjecture, holomorphy conjecture, and embedded Nash problem for Pfaffian ideals}

\begin{document}
	
	\maketitle
	
	\begin{abstract}
		We resolve the monodromy conjecture, holomorphy conjecture, and embedded Nash problem for Pfaffian ideals.
	\end{abstract}
    
    \tableofcontents
    
	\section{Introduction}
	
	Let $X$ be a $d$-dimensional smooth irreducible complex variety and $\mathcal I \subset \mathcal O_X$ be an ideal with a non-empty zero locus. The \textit{monodromy conjecture} predicts that exponents of poles of the motivic zeta function $Z_{\mathcal I}^{\mathrm{mot}}(T)$ are Verdier monodromy eigenvalues of $\mathcal I$, see \cite{PV10} or Section \ref{Preliminaries}. It was formulated by Denef and Loeser in the topological and motivic settings, see \cite{DL92,DL98}, and by Igusa in the $p$-adic setting, see \cite{Igu}, initially for a principal ideal $\mathcal{I} = (f)$. By linking the motivic and analytic aspects of singularities of $\mathcal{I}$, monodromy conjecture establishes its central role in singularity theory.
    

    Since its proposal decades ago, the monodromy conjecture has been the subject of extensive research. Significant progress has been made in various special cases: For a principal ideal $\mathcal{I}$, the $d=2$ case was addressed by the works of Loeser, Veys, and Blanco, see \cite{Loe88,Vey90,Bla24}. Esterov, Lemahieu, and Takeuchi solved the conjecture in the topological setting for Newton non-degenerate polynomials when $d \leq 4$, see \cite{E+}. In the case $d=3$, Veys provided key insights into the survival conditions for poles of motivic zeta functions, see \cite{Vey06}. For general $\mathcal{I}$, Van Proeyen and Veys settled the topological case for $d=2$, see \cite{PV10}. Shi and Zuo established a multiplicative Thom-Sebastiani type result for the monodromy conjecture in \cite{SZ24}. Chen and Zuo proved the conjecture for determinantal ideals, based on an in-depth analysis of monodromy zeta functions, see \cite{CZ25}. Despite these advances, the conjecture remains widely open in its general form. For a comprehensive overview, we refer to \cite{Vey25}.


    Denef proposed an analog, the \textit{holomorphy conjecture}, suggesting that the generalized topological zeta function $Z_{f}^{\mathrm{top},(d_{0})}(s)$ is holomorphic if $d_{0}$ does not divide orders of any eigenvalues of monodromy actions on cohomologies of Milnor fibers, see \cite{Den}. Lemahieu and Van Proeyen generalized this conjecture to the ideal case for $Z_{\mathcal I}^{\mathrm{top},(d_{0})}(s)$ by replacing Milnor monodromy eigenvalues with Verdier monodromy eigenvalues, see \cite{VPL11}.

    The holomorphy conjecture seems equally difficult. A complete answer is known for the case $d=2$ in \cite{Vey93,VPL11}. In the case $d=3$, this conjecture is solved for $\mathcal I$ given by a homogeneous polynomial by Rodrigues and Veys \cite{RV01}, and for $\mathcal I$ given by a Newton non-degenerate polynomial by Castryck, Ibadula, and Lemahieu \cite{CIL17}. However, this problem is again open in general cases. 
	
	Intrinsically, the motivic zeta function
	\begin{displaymath}
		Z_{\mathcal I}^{\mathrm{mot}}(T) := \sum_{p \geq 0} [\mathscr X_p^p(X,\mathcal I)]\mathbb L^{-pd} T^p
	\end{displaymath}
	is a normalized generating series for the classes in the Grothendieck ring of complex varieties of the contact loci, where $\mathbb L = [\mathbb A^1]$, see Section \ref{Preliminaries}. 
	Hence, a detailed understanding of contact loci is essential. To this end, the \textit{embedded Nash problem}, proposed in \cite[Section 2]{ELM04}, seeks a geometrical characterization of the irreducible components of the contact loci. It is a complementary question to the classical Nash problem, see \cite{IK03,FD16}. A complete solution to the embedded Nash problem is known for ideals of plane curves \cite{BBPP24}, affine toric varieties given by a cone together with a toric invariant ideal \cite{Ish04}, determinantal ideals  \cite{Roi}, and hyperplane arrangements \cite{BT20}, whereas the general cases have only been partially addressed through the application of dlt valuations, see \cite{BBPP24}.

	
    This paper focuses on the space $X = \mathcal{M} := M_m^{\mathrm{skew}}(\mathbb{C})$ of $m \times m$ complex skew-symmetric matrices, taking $\mathcal{I} = \mathcal{P}_{m_0}$ to be the ideal of $m_0 \times m_0$ Pfaffians of the $m\times m$ generic skew-symmetric matrix, with $n = \lfloor m/2 \rfloor$ and $1 \leq m_0 \leq n$. In literature, the monodromy conjecture for $(\mathcal M,\mathcal P_{m_0})$ is solved for the case $(m,m_0) = (2n+1,n)$, see \cite{LRWW17}, while the holomorphy conjecture and the embedded Nash problem remain completely open. In this work, we address all these three problems for general $(m,m_0)$.
    
	
	\subs{\bf Monodromy conjecture and holomorphy conjecture.}
	The first part of this paper relates the poles of motivic zeta functions, the holomorphy of generalized topological zeta functions, and the Verdier monodromy eigenvalues for $(\mathcal M,\mathcal P_{m_{0}})$.
	\begin{theorem}[Monodromy conjecture and holomorphy conjecture for Pfaffian ideals]\label{ThmA} 
    Let $\mathcal{M}=M_m^{\mathrm{skew}}(\mathbb C)$ be the space of $m\times m$ skew-symmetric matrices and $\mathcal P_{m_0}$ be the $m_{0}$-th Pfaffian ideal for every $1\le m_{0}\le \lfloor m/2 \rfloor$, then:
    
    \begin{enumerate}
        \item The set of poles of $Z_{\mathcal P_{m_0}}^{\mathrm{mot}}(T)$ is
		\begin{displaymath}
			\{-\frac{m(m-1)}{2m_0},-\frac{(m-2)(m-3)}{2(m_0-1)},\cdots, -\frac{(m-2m_0+2)(m-2m_0+1)}{2}\}.
		\end{displaymath}
        \item Let $d_{0}$ be a positive integer, then the poles of $Z_{\mathcal P_{m_0}}^{\mathrm{top},(d_{0})}(s)$ are contained in
		\begin{displaymath}
			\{-\frac{(m-2i+2)(m-2i+1)}{2(m_0-i+1)}\mid 1\leq i\leq m_0\text{ and }d_{0}\text{ divides } m_0-i+1\}.
		\end{displaymath}
        \item If $m_0 = \lfloor m/2\rfloor$, then $1$ is a Verdier monodromy eigenvalue of $\mathcal P_{m_0}$. Otherwise, for all $k = 1,\dots,m_0$ and $j = 0,1,\dots,m_0-k$, $e^{\frac{2j\pi \sqrt{-1}}{m_0-k+1}}$ is a Verdier monodromy eigenvalue of $\mathcal P_{m_0}$. 
    \end{enumerate}

		In particular, the monodromy conjecture and holomorphy conjecture hold for $(\mathcal M,\mathcal P_{m_0})$.
	\end{theorem}
    \begin{remark}
        If $m_0 = \lfloor m/2\rfloor$, then all poles of $Z_{\mathcal{P}_{m_0}}^{\mathrm{mot}}(T)$ are integers.
    \end{remark}

	\subs{\bf Embedded Nash problem.} Let $p \in \mathbb Z_{\geq 0}$ and $l \in \overline{\mathbb Z}_{\geq p} := \mathbb Z_{\geq p} \cup \{\infty\}$. Recall that $\mathscr X_p^l(\mathcal I)$ is the $p$-th contact loci at $l$-jet level, consisting of $l$-jets whose order of contact with $\mathcal I$ equals $p$. Here, we denote arcs by $\infty$-jets. The contact loci of $(\mathcal M,\mathcal P_{m_0})$ can be characterized combinatorially. Set
	\begin{align*}
		&\Lambda_{n,l} := \{(\lambda_1,\dots,\lambda_n) \in \overline{\mathbb Z}_{\geq 0}^{n} \mid \lambda_1 \leq \dots \leq \lambda_n \leq l+1\},\\
		&\Lambda_{n,l}^{m_0,p} := \{(\lambda_1,\dots,\lambda_n) \in \Lambda_{n,l} \mid \lambda_1+\dots+\lambda_{m_0} = p\}, \text{ and }\\
		& \Lambda_{n,l}^{m_0,p,*} := \{(\lambda_1,\dots,\lambda_n) \in \Lambda_{n,l}^{m_0,p} \mid 
		\lambda_{m_0-\lfloor p/k\rfloor+2} = \dots = \lambda_{n} = k, \ \lambda_{1} = \dots = \lambda_{m_0-\lfloor p/k\rfloor} = 0, \\
		&  \qquad\qquad \qquad \text{and}\ 
		\lambda_{m_0-\lfloor p/k\rfloor+1} = p- k\lfloor p/k\rfloor, \ k = \lceil p/m_0 \rceil,\dots,m_0 \}.
	\end{align*}
	Then for each $\bm \lambda = (\lambda_1,\dots,\lambda_n) \in \Lambda_{n,l}$, one can assign an irreducible Zariski locally closed subset $\mathcal C_{\bm \lambda,l}$ of $\mathcal M_l$, where we denote by $(-)_l$ the $l$-jet space of a variety. These subsets are orbits of the natural $(\mathrm{GL}_m(\mathbb C))_l$-action on $\mathcal M_l$, see Section \ref{contact_loci}.
	\begin{theorem}[Embedded Nash problem for Pfaffian ideals]\label{ThmB} Let $p \in \mathbb Z_{\geq 0}$ and $l \in \overline{\mathbb Z}_{\geq p}$, then:
	\begin{enumerate}
	    \item We have the following decomposition:
		\begin{displaymath}
			\mathscr X_p^l(\mathcal M,\mathcal P_{m_0}) = \bigsqcup_{\bm \lambda \in \Lambda_{n,l}^{m_0,p}} \mathcal C_{\bm \lambda,l}.
		\end{displaymath}
        \item For $\bm \lambda,\bm \lambda' \in \Lambda_{n,\infty}$, the closure of ${\mathcal C_{\bm \lambda,\infty}}$ contains $\mathcal C_{\bm \lambda',\infty}$ if and only if
		\begin{displaymath}
			\lambda_1+\dots+\lambda_k \leq \lambda_1'+\dots+\lambda_k'
		\end{displaymath} 
		for all $k = 1,\dots,n$.
	\end{enumerate}
    
		In particular, 
		\begin{displaymath}
			\mathscr X_p^{\infty}(\mathcal M,\mathcal P_{m_0}) = \bigcup_{\bm \lambda \in \Lambda_{n,\infty}^{m_0,p,*}} \overline{\mathcal C_{\bm \lambda,\infty}}
		\end{displaymath}
		is the decomposition into irreducible components. 
	\end{theorem}
    
    We hereby resolve the embedded Nash problem for $(\mathcal M,\mathcal P_{m_0})$.

	\subs{\bf Organization.} Section \ref{Preliminaries} is for preliminary materials on the contact loci, monodromy conjecture, and holomorphy conjecture. In Section \ref{contact_loci}, we determine the contact loci of Pfaffian ideals, solve the embedded Nash problem for them, and compute their topological zeta functions. In Section \ref{MC}, we prove the monodromy conjecture and the holomorphy conjecture for Pfaffian ideals.

	
		
	
	\section{Preliminaries}\label{Preliminaries}
	
	\subs{\bf Jet/arc spaces and contact loci.}\label{jet_and_arc} Let $X$ be a complex algebraic variety. For  $l\in \overline{\mathbb Z}_{\geq 0}$, we denote by $X_l$ the $l$-jet space of $X$, that is, the set of morphisms $\gamma : \mathrm{Spec}\, \mathbb C[[t]]/(t^{l+1}) \to X$ over $\mathbb C$. Each $X_l$ has the structure of a complex algebraic variety, but we only consider its reduced structure. When $l = \infty$, $X_{\infty}$ is called the arc space of $X$, consisting of arcs in $X$, that is, morphisms $\gamma : \mathrm{Spec}\, \mathbb C[[t]] \to X$ over $\mathbb{C}$, also endowed with the natural reduced $\mathbb C$-scheme structure. 
    We also denote arcs by $\infty$-jets. We have the natural truncation morphisms $\varphi_{l,l'}:X_l\to X_{l'}$ for $l\ge l'$ and $\varphi_{\infty, l}:X_{\infty}\to X_l$. A morphism of varieties $\mu:Y\to X$ induces morphisms $\mu_l:Y_l\to X_l$ compatible with the truncation maps. For $\gamma\in X_l$ with $l\in\overline{\mathbb Z}_{\geq 0}$, we denote by $\gamma(0):=\varphi_{l,0}(\gamma)\in X$ the   center of $\gamma$.
        
    Now suppose $X$ is a $d$-dimensional smooth complex variety and $\mathcal I$ is a coherent ideal sheaf of $X$ such that its zero locus $V(\mathcal I)$ is non-empty. Since $X$ is smooth, the truncations $\varphi_{l,l'}$ are Zariski locally trivial $\mathbb A^{d(l-l')}$-fibrations, see \cite[Proposition 3.7.5]{CLNS}. For $p \in \mathbb Z_{> 0}$ and $l \in \overline{\mathbb Z}_{\geq p}$, we define the $p$-th contact loci of $\mathcal I$ at the $l$-jet level:
        	\begin{displaymath}
        		\mathscr X_p^l(\mathcal I) := \{\gamma \in  X_p \mid \mathrm{ord}_\gamma (\mathcal I) = l\} = \varphi_{l,p-1}^{-1}(V(\mathcal I)_{p-1}) \setminus \varphi_{l,p}^{-1}(V(\mathcal I)_p),
        	\end{displaymath}
         where $\mathrm{ord}_{\gamma}(\mathcal I)$ is the order in $t$ of $\mathcal I$ via $\gamma$, and we set $\mathscr X_0^l(\mathcal I) =\varphi_{l,0}^{-1}(X\setminus V(\mathcal I))$. Then $\mathscr X^l_p(\mathcal I)=\varphi_{l,p}^{-1}(\mathscr X_p^p(\mathcal I))$ for all $l \in \overline{\mathbb Z}_{\geq p}$.

        Next, we recall a package of structure theorems of contact loci. Let $\mu : Y \to X$ be a log resolution of $\mathcal I$ such that $\mu : \mu^{-1}(X\setminus V(\mathcal I)) \to X\setminus V(\mathcal I)$ is an isomorphism. Suppose $\mathcal I\cdot \mathcal O_Y = \mathcal O_Y(-\sum_{i\in S} N_i E_i)$ and the canonical divisor $K_{\mu} = \sum_{i\in S} (\nu _i-1) E_i$, where $\bigcup_{i\in S} E_i = \mu^{-1}(V(\mathcal I))$ is a simple normal crossings divisor. Let $p\in\mathbb Z_{>0}$ and $l \in \overline{\mathbb Z}_{\geq p}$. We denote the $p$-th contact loci of $\mathcal{J}$ at $l$-jet level by $\mathscr{Y}_{\bm p}^l(\mathcal{J})$, where $\mathcal{J}$ is a coherent ideal sheaf of $Y$. For $I \subset S$ and $\bm p = (p_i)_{i\in S} \in \{0,1,\dots,l\}^{S}$, we denote $\mathscr Y_p^l(E_i) = \mathscr Y_p^l(\mathcal O_Y(-E_i))$, $\mathscr Y_{\bm p}^l(\bm E) = \bigcap_{i\in S} \mathscr Y_{p_i}^l(E_i)$, and $\mathscr X_{\bm p}^l(\mathcal I) = \mu_l(\mathscr Y_{\bm p}^l(\bm E))$.

        \begin{theorem}[\cite{ELM04,Cohomology_of_Contact_Loci}]\label{structure_thm_of_contact_loci}
            Notations are as above, then:
            \begin{enumerate}
                \item If $l' \geq l \geq \bm p \cdot \bm N = \sum_{i\in S} p_i N_i$, then $\varphi_{l',l}(\mathscr X^{l'}_{\bm p}(\mathcal I))=\mathscr X_{\bm p}^l(\mathcal I)$ and  $\mathscr X^{l'}_{\bm p}(\mathcal I) = \varphi_{l',l}^{-1}(\mathscr X_{\bm p}^l(\mathcal I))$.
                \item Each $\mathscr X^l_{\bm p}(\mathcal I)$ is smooth, irreducible, and locally closed in $X_l$. 
                \item There is a partition $\mathscr X_p^l(\mathcal I) = \bigsqcup_{\bm p\cdot \bm N = p} \mathscr X_{\bm p}^l(\mathcal I)$.
                \item The restricted morphism $\mu_l : \mathscr Y_{\bm p}^l(\bm E) \to \mathscr X_{\bm p}^l(\mathcal I)$ is a Zariski locally trivial $\mathbb A^{\bm p\cdot (\bm \nu-\bm{1}) }$-fibration, where $\bm p\cdot (\bm \nu-\bm{1}) = \sum_{i\in S} p_i(\nu_i-1)$.
                \item If $l \geq \bm p \cdot \bm N$, then the morphism $\mathscr Y_{\bm p}^l(\bm E) \to E_I^\circ$, $\tilde\gamma \mapsto \tilde\gamma(0)$, is a Zariski locally trivial $\mathbb A^{dl-\vert \bm p\vert} \times (\mathbb C^*)^{\vert I \vert}$-fibration, where $I = \{i\in S \mid p_i \neq 0\}$, $E_I^\circ = (\bigcap_{i\in I} E_i) \setminus (\bigcup_{j\in S\setminus I} E_j)$, and $\vert \bm p\vert = \sum_{i\in S} p_i$.
            \end{enumerate}
    	\end{theorem}

    \subs{\bf Grothendieck ring of complex varieties.}
    
	\begin{definition}[Grothendieck ring]\label{Ch2_Sec2_Grothendieck_Ring}
		Let $\mathrm{Var}_{\mathbb C}$ be the category of $\mathbb C$-varieties. The Grothendieck group of $\mathbb C$-varieties, denoted by $K_0(\mathrm{Var}_{\mathbb C})$, is defined to be the quotient group of the free abelian group with basis $\{[X]\}_{X\in\mathrm{Var}_{\mathbb C}}$, modulo the following relations:
		\begin{align*}
			& [X]-[Y], \textnormal{ if } X \simeq Y,\\
			& [X]-[X_{\mathrm{red}}],\\
			& [X]-[U]-[X\setminus U],\textnormal{ for any open set } U\subseteq X.
		\end{align*}
		One can further define a multiplication structure on $K_0(\mathrm{Var}_{\mathbb C})$ by
		\begin{displaymath}
			[X] \cdot [Y] := [X\times Y].
		\end{displaymath}
		This makes $K_0(\mathrm{Var}_{\mathbb C})$ a ring, called the Grothendieck ring of $\mathbb{C}$-varieties. Let $\mathbb L = [\mathbb A_{\mathbb C}^1]$ and $K_0(\mathrm{Var}_{\mathbb C})_{\mathbb L}$ be the localization at $\mathbb L$.
	\end{definition}

    The Grothendieck ring admits the virtual Poincar\'e specialization. It is the ring homomorphism $\mathrm{VP} : K_0(\mathrm{Var}_\mathbb C) \to \mathbb Z[w]$ defined by  $$[X] \to\sum_{m,i}  (-1)^{i+m} \dim \mathrm{Gr}_m^W H_c^i(X,\mathbb C) \cdot w^{m}$$ 
	for a variety $X$,
	where $W$ is the weight filtration. In particular, $\mathrm{VP}(\mathbb L) = w^2$. We will also denote by $\mathrm{VP}:K_0(\mathrm{Var}_\mathbb C)_{\mathbb L}\to \mathbb Z[w]_w$ the ring homomorphism induced by inverting $\mathbb L$ and $w$.

    One can further specialize $\mathbb Z[w]_w$ to $\mathbb Z$ by $\chi : w \mapsto -1$. Its composition with $\mathrm{VP}$, also denoted by $\chi$, is given by taking Euler characteristics, i.e. $\chi(X) = \sum_{i} (-1)^{i} \dim H_c^i(X,\mathbb C)$. 
    
	\subs{\bf Zeta functions and the monodromy conjecture.} We keep the notations in \ref{jet_and_arc}. The motivic zeta function of $\mathcal{I}$ can be computed via any log resolution $\mu$, cf. \cite{DL98,CLNS}:
    \begin{equation}\label{Zmot}
		Z_{\mathcal I}^{\mathrm{mot}}(T) = \sum_{I\subset S} [E_I^\circ] \cdot \prod_{i\in I} \frac{(\mathbb L-1)\mathbb L^{-\nu_i}T^{N_i}}{1-\mathbb L^{-\nu_i}T^{N_i}} \in K_0(\mathrm{Var}_\bC)_{\mathbb L}[[T]].
	\end{equation}
    Let $T = \mathbb L^{-s}$, then one defines the poles of $Z_{\mathcal I}^{\mathrm{mot}}$ in the expected way, see \cite{Vey25}. Applying the Euler characteristic $\chi$ to $Z_{\mathcal I}^{\mathrm{mot}}(T)$, then one specializes $Z_{\mathcal I}^{\mathrm{mot}}(T)$ to the topological zeta function
    \begin{equation}\label{Ztop}
        Z_{\mathcal I}^{\mathrm{top}}(s) = \sum_{I\subset S} \chi(E_I^\circ) \cdot \prod_{i\in I} \frac{1}{N_i s+\nu_i} \in \mathbb C(s).
    \end{equation}
    In particular, $Z_{\mathcal I}^{\mathrm{top}}(s)$ has fewer poles than $Z_{\mathcal I}^{\mathrm{mot}}(T)$.
    
    The other side of the monodromy conjecture is the set of Verdier monodromy eigenvalues $E({\mathcal I})$, cf. \cite{Ver83}, which are generalizations of monodromy eigenvalues of Milnor fibers, see \cite{BVWZ21,BSZ25} for an introduction. They can be computed using the Bernstein-Sato polynomial for ideals, see \cite{BMS06}. Or one can use the monodromy zeta function in the following way. Let $h : \widehat X \to X$ be the blow-up at $\mathcal I$. Since $\mathcal I \cdot \mathcal O_Y$ is an invertible sheaf, $\mu$ uniquely factors through $h$. We write the factorization as $\phi : Y \to \widehat{X}$.
    \begin{theorem}[\cite{PV10}] Let $\widehat{E} = h^{-1}(V(\mathcal I))$ and $e \in \widehat{E}$, then:
    \begin{enumerate}
        \item The monodromy zeta function of $\mathcal I$ at $e$ is
        \begin{displaymath}
            Z_{\mathcal I,e}^{\mathrm{mon}}(t) = \prod_{i\in S} (1-t^{N_i})^{\chi(E_i^\circ \cap \phi^{-1}(e))}.
        \end{displaymath}
        \item Denote by $\mathrm{PZ}(Z_{\mathcal I,e}^{\mathrm{mon}})$ the union of zeros and poles of $Z_{\mathcal I,e}^{\mathrm{mon}}$, then $E(\mathcal I) = \bigcup_{e\in \widehat{E}} \mathrm{PZ}(Z_{\mathcal I,e}^{\mathrm{mon}})$.
    \end{enumerate}
    
        In particular, all Verdier monodromy eigenvalues are roots of unity.
    \end{theorem}

    The monodromy conjecture for $\mathcal I$ is stated as follows.
    \begin{conjecture}[Monodromy conjecture, \cite{DL98,PV10}]
        Let $X$ be a smooth irreducible complex variety and $\mathcal I$ be a coherent ideal sheaf of $X$ with a non-empty zero locus. Suppose $s_0$ is a pole of $Z_{\mathcal I}^{\mathrm{mot}}(T)$, then $e^{2\pi \sqrt{-1}s_0} \in E(\mathcal I)$.
    \end{conjecture}

    \subs{\bf Generalized topological zeta functions and the holomorphy conjecture.}\label{gtzfhc} We keep the notations in \ref{jet_and_arc}. For a positive integer $d_{0}$, define the (generalized) topological zeta function of $\mathcal I$ with respect to $d_{0}$ as
    \begin{displaymath}
        Z^{\mathrm{top},(d_{0})}_{\mathcal I}(s) = \sum_{I\subset S_{d_{0}}} \chi(E_I^\circ) \cdot \prod_{i\in I} \frac{1}{N_i s+\nu_i} \in \mathbb C(s),
    \end{displaymath}
    where $S_{d_{0}} = \{i\in S \mid d_{0} \text{ divides } N_i\}$. By the weak factorization theorem \cite{AKW}, one verifies that $Z^{\mathrm{top},(d_{0})}_{\mathcal I}(s) $ is independent of the choice of log resolutions.

    In our notation, the holomorphy conjecture \cite[Conjecture 4]{VPL11} can be formulated as follows:
    \begin{conjecture}[Holomorphy conjecture]\label{HC}
        Let $X$ be a smooth irreducible complex variety and $\mathcal I$ be a coherent ideal sheaf of $X$ with a non-empty zero locus. Suppose $d_{0}$ is a positive integer that does not divide the order of any Verdier monodromy eigenvalues of $\mathcal I$, then  $Z^{\mathrm{top},(d_{0})}_{\mathcal I}(s)$ is holomorphic.
    \end{conjecture}
    
	
	\section{Contact loci} \label{contact_loci} Adopt the notations in the introduction. 
    Recall that $\mathcal M$ is the space of skew-symmetric $m\times m$ matrices with entries in $\mathbb C$. Let $R = \mathbb C[x_{ij}]/(x_{ij}+x_{ji} \mid 1\leq i,j\leq m)$ be the coordinate ring of $\mathcal M$, then $\mathcal P_{m_{0}}$ is the ideal of $R$ generated by Pfaffians of all $2m_{0}\times 2m_{0}$ principal submatrices of
	\begin{displaymath}
		\begin{pmatrix}
			0 & x_{12} & x_{13} & \cdots & x_{1,m-1} & x_{1m} \\
			-x_{12} & 0 & x_{23} & \cdots & x_{2,m-1} & x_{2m} \\
			 \vdots & \vdots & \vdots & \ddots & \vdots & \vdots \\
			 -x_{1m} & -x_{2m} & -x_{3m} & \cdots & -x_{m-1,m} & 0 
		\end{pmatrix}.
	\end{displaymath}


    \subs{\bf Classification of orbits.} We denote by $G$ the group $\mathrm{GL}_m(\mathbb C)$, which acts on $\mathcal M$ naturally via the congruence transformation:
	\begin{displaymath}
		g \cdot A := gAg^{T},\ g\in G, A \in \mathcal M.
	\end{displaymath}
	For $l \in \overline{\mathbb Z}_{\geq 0}$, $\mathcal M_l$ consists of skew-symmetric matrices of size $m\times m$ with entries in $\mathbb C[[t]]/(t^{l+1})$ and $G_{l}$ is the group given by invertible $m\times m$ matrices with entries in the same ring. The $G$-action on $\mathcal M$ also induces actions at the level of jets and arcs, again via congruence transformations. Motivated by \cite{Roi}, we start with the classification of $G_l$-orbits of $\mathcal M_l$. 
	\begin{definition}
		We associate with each $\bm \lambda := (\lambda_1,...,\lambda_n) \in \Lambda_{n,l}$ a matrix
		\begin{displaymath}
			\delta_{\bm \lambda,l} := \begin{pmatrix}
				t^{\lambda_1}J \\
				& t^{\lambda_2}J \\
				  && \ddots \\
				  &&& t^{\lambda_n}J \\
				  &&&& (0)
			\end{pmatrix} \in \mathcal M_{l},
		\end{displaymath}
		where $J = \begin{pmatrix}
			0 & 1 \\ -1 & 0
		\end{pmatrix}$, the $(0)$ is placed only when $m$ is odd, and we adopt the convention that $t^{\infty} = 0$. We denote the $G_{l}$-orbit of $\delta_{\bm \lambda,l}$ by $\mathcal C_{\bm \lambda,l}$. 
	\end{definition}
	\begin{proposition}\label{p3.1}
		Let $p\in \mathbb Z_{>0}$, $l \in \overline{\mathbb Z}_{\geq p}$, and $1\leq k \leq n$. Then every $G_{l}$-orbit of $\mathcal M_{l}$ is of the form $\mathcal C_{\bm \lambda,l}$ for some $\bm \lambda \in \Lambda_{n,l}$ and we have the following decomposition:
		\begin{displaymath}
			\mathscr{X}_{p}^{l}(\mathcal M, \mathcal P_k) = \bigsqcup_{\bm \lambda \in {\Lambda_{n,l}^{k,p}}} \mathcal C_{\bm \lambda,l}.
		\end{displaymath}
        In particular, Theorem \ref{ThmB}(1) holds.
	\end{proposition}
	\begin{proof}
		The first assertion follows from the facts that all matrices in $\mathcal M_{\infty}$ can be transformed to some $\delta_{\bm \lambda}$ under the action of $G_{\infty}$ and that the $k$-th Pfaffian ideal of $\mathcal O_{\mathcal M_{\infty}}$ stays invariant under congruence transformations.

        Therefore, it suffices to prove that $\mathcal{C}_{\bm{\lambda},l}=\mathcal{C}_{\bm{\lambda}',l}$ implies $\bm{\lambda}=\bm{\lambda}'$. Suppose $\bm \lambda \neq \bm \lambda'$ while $\mathcal{C}_{\bm{\lambda},l}=\mathcal{C}_{\bm{\lambda}',l}$, then there exists $A\in G_{l}$ such that $\delta_{\bm{\lambda},l}=A\delta_{\bm{\lambda}',l}A^{T}$. Let $u$ be the minimal number such that $\lambda_{u}\neq \lambda_{u}'$. Without loss of generality, we can suppose $\lambda_{u}<\lambda_{u}'$. Passing to the $\lambda_u$-jet level, we have
        \begin{displaymath}
            A
            \begin{pmatrix}
				t^{\lambda_1}J \\
				& \ddots \\
				  && t^{\lambda_{u-1}}J \\
                &&& 0\cdot J \\
                &&&& O
			\end{pmatrix}A^{T}=
            \begin{pmatrix}
				t^{\lambda_1}J \\
				& \ddots \\
				  && t^{\lambda_{u-1}}J \\
                &&& t^{\lambda_{u}}J  \\
                &&&& O
			\end{pmatrix},
        \end{displaymath}
        where $O$ is the zero matrix of size $(m-2u)\times (m-2u)$. Suppose $A = \begin{pmatrix}
            B & * \\
            * & *
        \end{pmatrix}$, where $B \in {M}_{2u}(\mathbb C[[t]]/(t^{\lambda_u+1}))$, then we have
        \begin{displaymath}
            B
            \begin{pmatrix}
				t^{\lambda_1}J \\
				& \ddots \\
				  && t^{\lambda_{u-1}}J \\
                &&& 0\cdot J
			\end{pmatrix}B^{T}=
            \begin{pmatrix}
				t^{\lambda_1}J \\
				& \ddots \\
				  && t^{\lambda_{u-1}}J \\
                &&& t^{\lambda_{u}}J  
			\end{pmatrix}.
        \end{displaymath}
        We write $B$ into the form $B=
        \begin{pmatrix}
            C & C' \\
            C'' & D
        \end{pmatrix}$ where $C$ is of size $(2u-2)\times (2u-2)$ and denote $S = \begin{pmatrix}
            t^{\lambda_1}J \\
				& \ddots \\
				  && t^{\lambda_{u-1}}J \\
        \end{pmatrix}$. Then 
        \begin{displaymath}
            B \begin{pmatrix}
                S &\\
                & 0\cdot J 
            \end{pmatrix}
            B^T = \begin{pmatrix}
                CSC^T & CS{C''}^T\\
                C'' S C^T & C'' S {C''}^T
            \end{pmatrix} = \begin{pmatrix}
                S & \\
                & t^{\lambda_u}J
            \end{pmatrix}.
        \end{displaymath}
        It follows that $CSC^T = S$ and $CS{C''}^T = 0$. Using Lemma \ref{solASAT=S1}(1), we obtain that $C$ is invertible. Then $CS{C''}^T = 0$ implies $S{C''}^T = 0$. Therefore, $0 = C''S{C''}^T \neq t^{\lambda_u}J$, a contradiction.
	\end{proof}

	\subs{\bf Jets in the stabilizer of $\mathcal C_{\bm \lambda,l}$.} \label{jets/arcs_in_the_stablizer}	Fix $m_{0},p \in \mathbb Z_{> 0}$, $l \in {\mathbb Z}_{\geq p}$, and $\bm \lambda \in \Lambda_{n,l}^{m_{0},p}$. Let $H_{\bm \lambda,l}$ be the stabilizer of $\delta_{\bm \lambda,l}$ in $G_{l}$.
	  Suppose $\lambda_1 = \dots = \lambda_{n_1} < \lambda_{n_1+1} = \dots =\lambda_{n_1+n_2} < \dots =\lambda_{n_1+\dots+n_q} < \lambda_{n_{1}+\dots+n_q+1} = l+1$. Here, we use the convention that $\lambda_{n+1} = \infty$. For $i=1,\dots,q$, let $N_i = \sum_{j=1}^i n_j$ and $\bar{\lambda}_{i}=\lambda_{N_{i}}$, then 
	\begin{displaymath}
		\delta_{\bm \lambda} = \begin{pmatrix}
			t^{\bar{\lambda}_{1}} S_1 \\
			& t^{\bar{\lambda}_{2}} S_2 \\
			&& \ddots\\
			&&& t^{\bar{\lambda}_{q}}S_q \\
			&&&& O
		\end{pmatrix},
	\end{displaymath}
	where $S_i$ is the block diagonal matrix given by $n_i$ copies of $J$ and $O$ is the zero matrix of size $(m-2N_q)\times (m-2N_q)$.

	\begin{lemma}\label{solASAT=S1}
		Let $\bar l\in \overline{\mathbb Z}_{\ge 0}$, $U_1,...,U_w$ be invertible matrices with entries in $\mathbb C[[t]]/(t^{\bar l+1})$ and $0 \leq \bar{\lambda}_{1},...,\bar{\lambda}_{w} \leq \bar l$ be a strictly increasing sequence of integers. Suppose $A$ is a $\mathbb C[[t]]/(t^{\bar l+1})$-valued square matrix such that 
		\begin{equation}\label{ASAT=S}
			A \begin{pmatrix}
				t^{{\bar \lambda}_1} U_1 \\
				& \ddots\\
				&& t^{{\bar \lambda}_w}U_w
			\end{pmatrix} A^T = 
			\begin{pmatrix}
				t^{{\bar \lambda}_1} U_1 \\
				& \ddots\\
				&& t^{{\bar \lambda}_w}U_w
			\end{pmatrix}.
		\end{equation}
	Let $r_i$ be the size of $U_i$ and represent $A$ as the block matrix $(A_{ij})_{1 \leq i,j \leq w}$, where each block $A_{ij}$ has size $r_i \times r_j$. Then:
    \begin{enumerate}
        \item $A,A_{11},\dots, A_{ww}$ are invertible. Equivalently, the constant terms of $\det A,\det A_{11},\dots, \det A_{ww}$ are all non-zero. Moreover, we have $A_{ii} U_i A_{ii}^T \equiv U_{i}$ modulo $(t)$. 
        \item  For any $i > j$, all entries of $A_{ij}$ are divisible by $t^{{\bar \lambda}_i-{\bar \lambda}_j}$.
    \end{enumerate}

	\end{lemma}
	\begin{proof}
		We proceed by induction on $w$. If $w = 1$, then (2) says nothing. Cancelling the $t^{\bar \lambda_1}$ on both sides of \eqref{ASAT=S}, we have $A_{11} U_1 A_{11}^{T} \equiv U_1$ modulo $(t^{\bar l+1-\bar{\lambda}_{1}})$. Note that $(\det A)^2 \cdot (\det U_1) = \det U_1 \in \mathbb C[[t]]/(t^{\bar l+1-{\bar \lambda}_1})$ and hence the constant term of $\det A$ is non-zero. This proves (1).
		
		Suppose the assertion holds for $w-1$. Let $A' = (A_{ij})_{1\leq i,j\leq w-1}$ and $U' = \begin{pmatrix}
			t^{{\bar \lambda}_1} U_1 \\
			& \ddots\\
			&& t^{{\bar \lambda}_{w-1}}U_{w-1}
		\end{pmatrix}$. Then $A' U' A'^T \equiv U'$ modulo $(t^{\bar{\lambda}_{w}})$. By the induction hypothesis, we have $A',A_{11},\dots,A_{w-1,w-1}$ are all invertible, $A_{ii}U_iA_{ii}^T \equiv U_{i}$ modulo $(t)$, and all entries of $A_{ij}$ are divisible by $t^{\bar \lambda_i-\bar \lambda_j}$ for $1\leq j< i\leq w-1$. Moreover, from \eqref{ASAT=S}, we have
		\begin{displaymath}
			A'U' \cdot \begin{pmatrix}
				t^{{\bar \lambda}_1}A_{w1}^T \\
				\vdots\\
				t^{{\bar \lambda}_w}A_{w,w-1}^T
			\end{pmatrix} = 0
		\end{displaymath}
		in $M_{(r_1+\dots+r_{w-1})\times r_w}(\mathbb C[[t]]/(t^{{\bar \lambda}_w}))$. Then (2) follows immediately since $A'$ and $U'$ are invertible. For (1), by comparing the last $r_w\times r_w$ block of the two sides of \eqref{ASAT=S}, we observe that
		\begin{displaymath}
			t^{{\bar \lambda}_w} A_{ww} U_w A_{ww}^T + t^{{\bar \lambda}_w}\sum_{i=1}^{w-1} t^{{\bar \lambda}_w-{\bar \lambda}_i} \widetilde A_{wi} U_i \widetilde A_{wi}^T = t^{{\bar \lambda}_w} U_w
		\end{displaymath}	
		in $M_{r_w}(\mathbb C[[t]]/(t^{\bar l+1}))$, where $\widetilde A_{wi} = t^{-({\bar \lambda}_w-{\bar \lambda}_i)}A_{wi} \in M_{r_w \times r_i}(\mathbb C[[t]]/(t^{\bar l+1-{\bar \lambda}_w+{\bar \lambda}_i}))$ and it is well defined by (2). It follows that $A_{ww}U_w A_{ww}^T \equiv U_w$ modulo $(t)$ and the constant term of $\det A_{ww}$ is thus non-zero. Now, since for any $i > j$, all entries of $A_{ij}$ are divisible by $t$, we have $\det A \equiv \det A' \cdot \det A_{ww}$ modulo $(t)$. This completes the proof of (1).
	\end{proof}
	
	
	\begin{lemma}\label{solASAT=S2}
		Notations are as in Lemma \ref{solASAT=S1}. Further assume $w = q$, $U_i = S_i$, $r_i = 2n_i$ for all $i$, and that $\bar l = l$ is an integer. Then the closed subvariety of $G_l$ of all $A$ satisfying \eqref{ASAT=S} is isomorphic to 
		\begin{displaymath}
			\mathbb A^{(2N_{q}^{2}-2\sum_{i=1}^{q} n_{i}^{2})(l+1)} \times \prod_{i < j} \mathbb A^{4{\bar \lambda}_{i}n_{i}n_{j}} \times \prod_{i=1}^{q} \mathrm{Sp}_{2n_{i}} \times \prod_{i=1}^{q} \mathbb A^{(l-{\bar \lambda}_{i})n_{i}(2n_i+1)} \times \prod_{i=1}^{q} \mathbb A^{4{\bar \lambda}_i n_{i}^{2}},
		\end{displaymath}
		where $\mathrm{Sp}_{2n_{i}}$ is the symplectic group of order $2n_{i}$.
	\end{lemma}
	\begin{proof}
		Let $A_{ij} = \sum_{u=0}^l A_{ij}^{(u)}t^u$ with $A_{ij}^{(u)} \in {M}_{2n_i}(\mathbb C)$. We will show (1) there is no constraint on $A_{ij}$ for $i< j$; (2) for any $i>j,0 \leq u \leq l-{\bar \lambda}_j$, $A_{ij}^{(u)}$ with is determined by $A_{i'j'}, A_{ik'}$, and $A_{ii}^{(0)},\dots, A_{ii}^{(u-{\bar \lambda}_i+{\bar \lambda}_j)}$, where $i' < i,j' = 1,\dots,w$, and $k' \geq i$; (3) there is no constraint on $A_{ij}^{(u)}$ for $i>j$ and $u \geq l-{\bar \lambda}_j+1$; (4) $A_{ii}^{(0)} \in \mathrm{Sp}_{2n_i}$; (5) for $1\leq u\leq l-{\bar \lambda}_i$, $A_{ii}^{(u)}$ is constrained by $A_{i'j'}$, $A_{ii}^{(0)},...,A_{ii}^{(u-1)}$, and $n_i(2n_i-1)$ independent linear equations, where $i'\leq i$ and $(i',j') \neq (i,i)$; (6) there is no constraint on $A_{ii}^{(u)}$ for $u \geq l-{\bar \lambda}_i + 1$. The conclusion follows by combining (1)-(6).

%
%
		
		We treat the entries of $A_{ij}$ as parameters for $i < j$ and as variables for $i \ge j$, and proceed to solve equation \eqref{ASAT=S}. Then it suffices to show (2)-(6) for fixed $A_{ij}$ with $i<j$.
		
		Let $\Lambda_i := \begin{pmatrix}
			A_{11} & \cdots & A_{1,i-1} \\
			\vdots & \ddots & \vdots\\
			A_{i-1,1} & \cdots & A_{i-1,i-1}  
		\end{pmatrix}$ and $\Theta_{i} := \begin{pmatrix}
		A_{1i} & \cdots & A_{1w} \\
		\vdots & \ddots & \vdots\\
		A_{i-1,i} & \cdots & A_{i-1,w}  
	\end{pmatrix}$. Then each $\Lambda_i$ is invertible by Lemma \ref{solASAT=S1}(1). For $i > j$, define $\widetilde A_{ij} := t^{{\bar \lambda}_j-{\bar \lambda}_i} A_{ij} \in M_{r_i \times r_j}(\mathbb C[[t]]/(t^{l+1-{\bar \lambda}_i+{\bar \lambda}_j}))$. By comparing the $(1,i)$-th, $(2,i)$-th, $\dots$, $(i-1,i)$-th blocks on both sides of \eqref{ASAT=S}, we obtain 
	\begin{equation}\label{r1r2ri-1}
			t^{{\bar \lambda}_i} \Lambda_i \cdot 
			\begin{pmatrix}
				S_1 \\
				& \ddots\\
				&& S_{i-1}
			\end{pmatrix}
			\cdot 
			\begin{pmatrix}
				\widetilde A_{i1}^T \\
				\vdots\\
				\widetilde A_{i,i-1}^T
			\end{pmatrix}
			=
			-t^{{\bar \lambda}_i}\Theta_i \cdot 
			\begin{pmatrix}
				S_i \\
				& \ddots\\
				&& t^{{\bar \lambda}_w-{\bar \lambda}_i}S_{w}
			\end{pmatrix}
			\cdot
			\begin{pmatrix}
				A_{ii}^T \\
				\vdots\\
				A_{iw}^T
			\end{pmatrix}
	\end{equation}
	in $M_{(r_1+\dots+r_{i-1})\times r_i}(\mathbb C[[t]]/(t^{l+1}))$. Since $\Lambda_i$ and $\begin{pmatrix}
		S_1 \\
		& \ddots\\
		&& S_{i-1}
	\end{pmatrix}$ are invertible, by comparing the order of $t$ on both sides of \eqref{r1r2ri-1}, we see that (2)(3) holds. 

	To prove (4)-(6), comparing the $(i,i)$-th block on the both sides of \eqref{ASAT=S}, one observes
	\begin{equation}\label{riri}
		\sum_{j=1}^{i-1} t^{{\bar \lambda}_i-{\bar \lambda}_j} \widetilde A_{ij} S_j \widetilde A_{ij}^T + A_{ii} S_i A_{ii}^T + \sum_{j=i+1}^w t^{{\bar \lambda}_j-{\bar \lambda}_i} A_{ij}S_j = S_i
	\end{equation}
	in $M_{r_i\times r_i}(\mathbb C[[t]]/(t^{l+1-{\bar \lambda}_i}))$. Since \eqref{riri} involves only terms in $t$ of degree less than $l-{\bar \lambda}_i$, this establishes (6). Reducing \eqref{riri} modulo $(t)$ implies that $A_{ii}^{(0)}$ lies in $\mathrm{Sp}_{2n_i}$, with no further constraints. For (5), we solve $A_{ii}^{(u)}$ inductively. Suppose $u \leq l-{\bar \lambda}_i$ and $A_{ii}^{(u-1)},\dots,A_{ii}^{(0)}$ are fixed, then $A_{{ii}}^{(u)}$ is given by
	\begin{equation}\label{ririu}
		A_{ii}^{(u)}S_i (A_{ii}^{(0)})^T + A_{ii}^{(0)} S_i (A_{ii}^{(u)})^T = t^{-u}\bigg( S_i - \sum_{j=1}^{i-1} t^{{\bar \lambda}_i-{\bar \lambda}_j} \widetilde A_{ij} S_j \widetilde A_{ij}^T - \sum_{j=i+1}^w t^{{\bar \lambda}_j-{\bar \lambda}_i} A_{ij}S_j\bigg) 
	\end{equation}
	in $M_{r_i\times r_i}(\mathbb C[[t]]/(t))$. We denote $W :=A_{ii}^{(u)}S_i (A_{ii}^{(0)})^T + A_{ii}^{(0)} S_i (A_{ii}^{(u)})^T$. By (2), $W$ is a matrix function in $A_{ij},i < j$ and $A_{ii}^{(0)},\dots,A_{ii}^{(u-1)}$. Then 
	\begin{equation}\label{ririu2}
			(A_{ii}^{(0)})^{-1}A_{ii}^{(u)}S_i + S_i \big((A_{ii}^{(0)})^{-1}A_{ii}^{(u)}\big)^{T} = (A_{ii}^{(0)})^{-1} W \big((A_{ii}^{(0)})^{-1}\big)^T  
	\end{equation}
	in $M_{2n_i\times 2n_i}(\mathbb C[[t]]/(t)) = M_{2n_i \times 2n_i}(\mathbb C)$. Suppose the matrix $(A_{ii}^{(0)})^{-1} W \big((A_{ii}^{(0)})^{-1}\big)^T$ is $(v_{i''j''})$. Let $B = (\beta_{i''j''}) := (A_{ii}^{(0)})^{-1}A_{ii}^{(u)}$, then $BS_i + S_iB^T = (c_{i''j''})$, where
	\begin{displaymath}
		c_{i''j''} = \begin{cases}
			\beta_{2a+1,2b+2}-\beta_{2b+1,2a+2}, & i'' = 2a+1, j'' = 2b+1,\\
			-\beta_{2a+1,2b-1}-\beta_{2b,2a+2}, & i'' = 2a+1, j'' = 2b,\\
			\beta_{2a,2b+2}+\beta_{2b+1,2a-1}, & i'' = 2a, j'' = 2b+1,\\
			-\beta_{2a,2b-1} + \beta_{2b,2a-1}, & i'' = 2a, j'' = 2b.
		\end{cases}
	\end{displaymath}
	Therefore, $\{\beta_{i''j''}\}$ is determined by $n_i(2n_i-1)$ linearly independent equations $c_{i''j''} = v_{i''j''}, i''<j''$, implying (5). 
	\end{proof}
	
	\begin{lemma}\label{solASAT=S3}
    Notations are as in the beginning of \ref{jets/arcs_in_the_stablizer} and we assume that $l$ is an integer. Then $H_{\bm \lambda,l}$, the stabilizer of $\delta_{\bm \lambda,l}$, is isomorphic to
		\begin{displaymath}
			H'_{\bm \lambda,l} \times \mathbb A^{2(l+1)(m-2N_q)N_q} \times \mathrm{GL}_{m-2N_q} \times \mathbb A^{l(m-2N_q)^2} \times \prod_{i=1}^q \mathbb A^{\bar{\lambda}_{i} n_i (m-2N_q)},
		\end{displaymath}
		where $H'_{\bm \lambda,l}$ is the stabilizer of 
		\begin{displaymath}
			\delta_{\bm \lambda,l}' := \begin{pmatrix}
				t^{\bar{\lambda}_{1}} S_1 \\
				& t^{\bar{\lambda}_{2}} S_2 \\
				&& \ddots\\
				&&& t^{\bar{\lambda}_{q}}S_q 
			\end{pmatrix}
		\end{displaymath}
		under the action of $(\mathrm{GL}_{2N_q}(\mathbb C))_l$.
	\end{lemma}
	\begin{proof}
		Let $A = \begin{pmatrix}
			A_1 & A_2 \\
			A_3 & A_4
		\end{pmatrix} \in H_{\bm \lambda,l}$, where $A_1,A_2,A_3$, and $A_4$ have size $2N_q\times 2N_q, 2N_q\times (m-2N_q), (m-2N_q)\times 2N_q$, and $(m-2N_q)\times (m-2N_q)$ respectively. Then
		\begin{equation}\label{deltagamm0}
			A\delta_{\bm \lambda,l}A^T = 
			\begin{pmatrix}
				A_1 \delta_{{\bm \lambda},l}' A_1^T & A_1\delta_{{\bm \lambda},l}' A_3^T \\
				A_3\delta_{{\bm \lambda},l}' A_1^T & A_3\delta_{{\bm \lambda},l}' A_3^T
			\end{pmatrix} = 
			\begin{pmatrix}
				\delta_{{\bm \lambda},l}' \\
				& O
			\end{pmatrix}.
		\end{equation}
		Therefore, $A_1\delta_{{\bm \lambda},l}' A_1^T \equiv \delta_{{\bm \lambda},l}'$ modulo $(t^{l+1})$. By Lemma \ref{solASAT=S1}, $A_1$ is invertible. Then we have $\delta_{{\bm \lambda},l}' A_3^T = 0$ and there is no other constraint from \eqref{deltagamm0}. Suppose $A_3 = (A_{31},\dots,A_{3q})$, where $A_{3i} \in M_{(m-2N_{q}) \times 2n_{i}}(\mathbb C[[t]]/(t^{l+1}))$. Then $\delta_{{\bm \lambda},l}' A_3^T = 0$ is equivalent to all entries of $A_{3i}$ being divisible by $t^{l+1-\bar{\lambda}_{i}}$. It follows that the constant term of $\det A$ equals the constant term of $\det A_1 \cdot \det A_4$. So $A_4$ is invertible. We thus conclude the proof.
	\end{proof}

    \subs{\bf Virtual Poincar\'e and topological zeta functions.} In this subsection, we will simplify the virtual Poincar\'e zeta function of $\mathcal P_{m_0}$ (i.e. the virtual Poincar\'{e} specialization of $Z_{\mathcal P_{m_0}}^{\mathrm{mot}}(T)$) and compute the topological zeta function. Let 
    $$Z_{\mathcal{P}_{m_{0}}}^{\mathrm{VP}}(T)=\sum_{p\geq 0}\mathrm{VP}(\mathscr{X}_{p}^{p}(\mathcal{M},\mathcal{P}_{m_{0}}))w^{-pm(m-1)}T^{p} \in \mathbb Z[w]_w[[T]] \subset \mathbb Z((w))[[T]]$$ 
    be the virtual Poincar\'e zeta function.

    By Proposition \ref{p3.1}, for each $l \in \mathbb Z_{\geq p}$, we have $$\mathrm{VP}(\mathscr{X}_{p}^{p}(\mathcal{M},\mathcal{P}_{m_{0}}))w^{-pm(m-1)} = \mathrm{VP}(\mathscr{X}_{p}^{l}(\mathcal{M},\mathcal{P}_{m_{0}}))w^{-lm(m-1)} = \sum_{\bm \lambda \in \Lambda_{n,l}^{m_0,p}} \mathrm{VP}(\mathcal C_{\bm \lambda,l})w^{-lm(m-1)}.$$
    With some adjustments, the formula also holds when $l = \infty$. Let $\widehat{\Lambda}_{n,p}^{m_0} = \Lambda_{n,p}^{m_0,\infty} \cap \mathbb Z_{\geq 0}^{m_0}$, and for $\bm \lambda \in \widehat{\Lambda}_{n,p}^{m_0}$ and $l > \lambda_1+\dots+\lambda_{n}$, set $\mu(\mathcal C_{\bm \lambda,\infty}) := \mathrm{VP}(\mathcal C_{\bm \lambda,l}) \cdot w^{-lm(m-1)}$. Then
    \begin{equation}\label{motivic_integration}
        \mathrm{VP}(\mathscr{X}_{p}^{p}(\mathcal{M},\mathcal{P}_{m_{0}}))w^{-pm(m-1)} = \sum_{\bm \lambda \in \widehat{\Lambda}_{n,p}^{m_0}} \mu(\mathcal C_{\bm \lambda,\infty}) \in \mathbb Z((w))[[T]],
    \end{equation}
    where by abuse of notation, we let $\bm{\lambda}$ also denote its truncation in $(\mathbb Z/(l+1)\mathbb Z)^{n}$. The equality can be directly checked using Proposition \ref{p3.1} and the computation below, or one may turn to the motivic integration, cf. \cite{CLNS}.

   Keep the notations in \ref{jets/arcs_in_the_stablizer}. For $\bm \lambda \in \widehat{\Lambda}_{n,p}^{m_0}$, we have a sequence of integers $\bm n(\bm \lambda) = (n_1,\dots,n_q)$ such that $\lambda_1 = \dots = \lambda_{n_1} < \lambda_{n_1+1} = \dots =\lambda_{n_1+n_2} < \dots =\lambda_{n_1+\dots+n_q} < \lambda_{n_{1}+\dots+n_q+1}$ and $n_1+\dots+n_q = n$. Recall the notations $N_i := \sum_{j=1}^i n_j$ and $\bar \lambda_i := \lambda_{N_i}$. Then by Lemma \ref{solASAT=S2} and Lemma \ref{solASAT=S3}, we have 
    \begin{align*}
        & [\mathcal{C}_{\bm{\lambda},l}]\mathbb{L}^{-\frac{lm(m-1)}{2}}\\
        & = [G_{l}/H_{\bm{\lambda},l}]\mathbb{L}^{-\frac{lm(m-1)}{2}}\\
        & =  \frac{[\mathrm{GL}_{m}]\mathbb{L}^{lm^{2}}\mathbb{L}^{-\frac{lm(m-1)}{2}}}{[H_{\bm{\lambda},l}'][\mathrm{GL}_{m-2n}]\mathbb{L}^{lm(m-2n)}\mathbb{L}^{2(m-2n)(n+\sum_{i=1}^{q}\bar{\lambda}_{i}n_{i})}}\\
        & =  \frac{[\mathrm{GL}_{m}]\mathbb{L}^{lm^{2}}\mathbb{L}^{-\frac{lm(m-1)}{2}}\mathbb{L}^{-lm(m-2n)}}{[\mathrm{GL}_{m-2n}]\prod_{i=1}^{q}[\mathrm{Sp}_{2n_{i}}]\mathbb{L}^{l(2n^{2}+n)}\mathbb{L}^{2n^{2}-2\sum_{i=1}^{q}n_{i}^{2}+\sum_{i=1}^{q}\bar{\lambda}_{i}n_{i}(4n-4N_{i}-1+2n_{i})}\mathbb{L}^{2(m-2n)(n+\sum_{i=1}^{q}\bar{\lambda}_{i}n_{i})}}\\
        & =  \frac{[\mathrm{GL}_{m}]}{[\mathrm{GL}_{m-2n}]\prod_{i=1}^{q}[\mathrm{Sp}_{2n_{i}}]\mathbb{L}^{2mn-2n^{2}-2\sum_{i=1}^{q}n_{i}^{2}}}\cdot \frac{\mathbb{L}^{-\frac{1}{2}l(m-2n)(m-2n-1)}}{\mathbb{L}^{\sum_{i=1}^{q}\bar{\lambda}_{i}n_{i}(2m-4N_{i}+2n_{i}-1)}}\\
        & =  \frac{[\mathrm{GL}_{m}]}{[\mathrm{GL}_{m-2n}]\prod_{i=1}^{q}[\mathrm{Sp}_{2n_{i}}]\mathbb{L}^{n-2\sum_{i=1}^{q}n_{i}^{2}}}\cdot \frac{1}{\mathbb{L}^{\sum_{i=1}^{q}(\bar{\lambda}_{i}+1)n_{i}(2m-4N_{i}+2n_{i}-1)}}.
    \end{align*}

    For $\bm{\lambda} \in \widehat{\Lambda}_{m_0,p}$, denote $b(\bm{\lambda})=\sum_{i=1}^{q}(\bar{\lambda}_{i}+1)n_{i}(2m-4N_{i}+2n_{i}-1)$, and for a tuple of positive integers $\bm u = (u_1,\dots,u_q)$ with $u_1+\dots+u_q = n$, let $Y_{\bm{u}}=\mathrm{VP}(\frac{[\mathrm{GL}_{m}]}{[\mathrm{GL}_{m-2n}]\prod_{i=1}^{q}[\mathrm{Sp}_{2u_{i}}]\mathbb{L}^{n-2\sum_{i=1}^{q}u_{i}^{2}}})$. Then the virtual Poincar\'{e} zeta function is 
    \begin{displaymath}
        Z_{\mathcal{P}_{m_{0}}}^{\mathrm{VP}}(T)=\sum_{p\geq 0}\sum_{\bm{\lambda}\in \widehat{\Lambda}_{n}^{m_{0},p}}\frac{Y_{\bm{n}(\bm{\lambda})}T^{\sum_{i=1}^{m_{0}}\lambda_{i}}}{w^{2b(\bm{\lambda})}}=\sum_{0\le \lambda_{1}\le \cdots \le \lambda_{n}<\infty}\frac{Y_{\bm{n}(\bm{\lambda})}T^{\sum_{i=1}^{m_{0}}\lambda_{i}}}{w^{2b(\bm{\lambda})}}
    \end{displaymath}
    and it can be further rewritten as
    \begin{displaymath}
        Z_{\mathcal{P}_{m_{0}}}^{\mathrm{VP}}(T)=\sum_{q=1}^{n}\sum_{u_{1}+\cdots +u_{q}= n}Y_{\bm u}\sum_{\bm{n}(\bm{\lambda})=(u_{1},\cdots ,u_{q})}\frac{T^{\sum_{i=1}^{m_{0}}\lambda_{i}}}{w^{2b(\bm{\lambda})}}.
    \end{displaymath}
    For a tuple $\bm{u}=(u_{1},\cdots ,u_{q})$, let $\alpha_{i}(\bm{u})=\max\{ 0,m_{0}-\sum_{j=1}^{i-1}u_{j} \}$ if $m_{0}<\sum_{j=1}^{i}u_{j}$, and $\alpha_{i}(\bm{u})=u_{i}$ otherwise. Then we have $\sum_{i=1}^{m_{0}}\lambda_{i}=\sum_{i=1}^{q}\alpha_{i}(\bm{n}(\bm{\lambda}))\bar{\lambda}_{i}$. Let $U_{i}=\sum_{j=1}^{i}u_{i}$. With the bijective map 
    \begin{displaymath}
        \tau:\mathbb{Z}_{>0}^{q} \to \{ (\bar{\lambda}_{1}+1,\cdots ,\bar{\lambda}_{q}+1)\in \mathbb{Z}_{>0}^{q}|\bar{\lambda}_{1}<\cdots <\bar{\lambda}_{q} \},(a_{1},\cdots ,a_{q})\mapsto (\sum_{i=1}^{1}a_{i},\cdots ,\sum_{i=1}^{q}a_{i}),
    \end{displaymath}
    we may write $\bar{\lambda}_{i}+1=\sum_{j=1}^{i}a_{j}$. Hence,
    \begin{align*}
        \sum_{\bm{n}(\bm{\lambda})= \bm{u}}\frac{T^{\sum_{i=1}^{m_{0}}\lambda_{i}}}{w^{2b(\bm{\lambda})}}
        & =\sum_{\bar{\lambda}_{1}<\cdots <\bar{\lambda}_{q}}\frac{T^{\sum_{i=1}^{q}\alpha_{i}(\bm{u})\bar{\lambda}_{i}}}{w^{2\sum_{i=1}^{q}(\bar{\lambda}_{i}+1)u_{i}(2m-4U_{i}+2u_{i}-1)}}\\
        & =\sum_{(a_{1},\cdots ,a_{q})\in \mathbb{Z}_{>0}^{q}}\frac{T^{-m_{0}+\sum_{i=1}^{q}(\sum_{j=i}^{q}\alpha_{j}(\bm{u}))a_{i}}}{w^{2\sum_{i=1}^{q}(\sum_{j=i}^{q}u_{j}(2m-4U_{j}+2u_{j}-1))a_{i}}}\\
        & =T^{-m_{0}}\prod_{i=1}^{q}(\sum_{a_{i}=1}^{+\infty}(\frac{T^{\sum_{j=i}^{q}\alpha_{j}(\bm{u})}}{w^{2\sum_{j=i}^{q}u_{j}(2m-4U_{j}+2u_{j}-1)}})^{a_{i}})\\
        & =T^{-m_{0}}\prod_{i=1}^{q}\frac{w^{-2\sum_{j=i}^{q}u_{j}(2m-4U_{j}+2u_{j}-1)}T^{\sum_{j=i}^{q}\alpha_{j}(\bm{u})}}{1-w^{-2\sum_{j=i}^{q}u_{j}(2m-4U_{j}+2u_{j}-1)}T^{\sum_{j=i}^{q}\alpha_{j}(\bm{u})}}.
    \end{align*}

    The simplification of $Z_{\mathcal{P}_{m_{0}}}^{\mathrm{VP}}(T)$ stops here. We now turn to the computation of the topological zeta function $Z_{\mathcal{P}_{m_{0}}}^{\mathrm{top}}(s)$. The key to the computation is the degree of $(w^2-1)$ in each term of $Z_{\mathcal{P}_{m_{0}}}^{\mathrm{VP}}(T)$. So we rewrite $Z_{\mathcal{P}_{m_{0}}}^{\mathrm{VP}}(T)$ as
    \begin{equation}\label{ZVP}
        \sum_{q=1}^{n}\sum_{u_{1}+\cdots +u_{q}=n}\frac{Y_{\bm u}T^{-m_{0}}}{(w^2-1)^{q}}\cdot \prod_{i=1}^{q}\frac{(w^2-1)w^{-2\sum_{j=i}^{q}u_{j}(2m-4U_{j}+2u_{j}-1)}T^{\sum_{j=i}^{q}\alpha_{j}(\bm{u})}}{1-w^{-2\sum_{j=i}^{q}u_{j}(2m-4U_{j}+2u_{j}-1)}T^{\sum_{j=i}^{q}\alpha_{j}(\bm{u})}}.
    \end{equation}
    Note that
    \begin{align*}
        [\mathrm{Sp}_{2k}] & = \mathbb L^{k^2}(\mathbb L^{2k}-1)(\mathbb L^{2k-2}-1)\cdots (\mathbb L^2-1),\\
        [\mathrm{GL}_k] & = \mathbb L^{\frac{k(k-1)}{2}}(\mathbb L^k-1) \cdots (\mathbb L-1).
    \end{align*}
    For $\bm{u}=(u_{1},\cdots ,u_{q})$ with $u_{1}+\cdots +u_{q}= n$, we have
    \begin{align*}
        Y_{\bm u} & =\frac{1}{w^{2n-4\sum_{i=1}^{q}u_{i}^{2}}}\cdot \frac{{w}^{m(m-1)}\prod_{i=1}^{m}(w^{2i}-1)}{(\prod_{i=1}^{m-2n}(w^{2i}-1))w^{2\sum_{i=1}^{q}u_{i}^{2}}\prod_{i=1}^{q}\prod_{j=1}^{u_{i}}(w^{4j}-1)}\\
        & =(w^2-1)^{m-(m-2n)-\sum_{i=1}^{q}u_{i}}\mathrm{W}_{\bm u}(w^2)\\
        & =(w^2-1)^{n}\mathrm{W}_{\bm u}(w^2),
    \end{align*}
    where $\mathrm{W}_{\bm u}(t) \in \mathbb C(t)$ is a rational function with no zeros or poles at $t=1$. Consequently, $\chi(Y_{\bm u}/(w^2-1)) = 0$ for all $(q,\bm u) \neq (n,\bm 1)$, where $\bm 1 := (1,\dots,1) \in \mathbb Z_{\geq 0}^n$. In this case, $\alpha_{i}(\bm 1)=0$ if $m_{0}<i$ and $\alpha_{i}(\bm 1)=1$ otherwise. Let $\epsilon_{i}=\mathrm{max}\{ 0,m_{0}-i+1 \}$, then the topological zeta function is
    \begin{align*}
         Z_{\mathcal{P}_{m_{0}}}^{\mathrm{top}}(s) & =\mathrm{W}_{\bm 1}(1)\prod_{i=1}^{n}\frac{1}{\sum_{j=i}^{n}\alpha_{j}(\bm 1)s+\sum_{j=i}^{n}(2m-4j+1)}\\
        & =\mathrm{W}_{\bm 1}(1)\prod_{i=1}^{n}\frac{1}{\epsilon_{i}s+(2m-2n-2i+1)(n-i+1)}\\
        & = \frac{m!}{2^{2n}}\prod_{i=1}^{n}\frac{1}{\epsilon_{i}s+(2m-2n-2i+1)(n-i+1)}.
    \end{align*}

    \begin{theorem}\label{poles_of_top_zeta}
        The set of poles of $Z_{\mathcal{P}_{m_0}}^{\mathrm{top}}(s)$ is
        \begin{displaymath}
    		\{ -\frac{m(m-1)}{2m_0},-\frac{(m-2)(m-3)}{2(m_0-1)},\cdots,-\frac{(m-2m_0+2)(m-2m_0+1)}{2} \}.
    	\end{displaymath}
    \end{theorem}
    
    \subs{\bf Embedded Nash problem.} In this subsection, we focus on the case $l = \infty$ and classify the irreducible components of $\mathscr{X}_p^{\infty}(\mathcal{P}_{m_0})$ for a fixed $p \in \mathbb Z_{>0}$. We define a partial order $\preceq$ in $\Lambda_{n,\infty}$ by 
    \begin{displaymath}
        \bm{\lambda} \preceq \bm{\lambda}' \iff \lambda_1 + \dots + \lambda_k \leq \lambda_1' + \dots + \lambda_k'\text{, for all }k. 
    \end{displaymath}
    \begin{proposition}
        For $\bm \lambda,\bm \lambda' \in \Lambda_{n,\infty}$, $\mathcal C_{\bm \lambda',\infty}$ is contained in the closure of $\mathcal C_{\bm \lambda,\infty}$ if and only if $\bm \lambda \preceq \bm \lambda'$. In particular, Theorem \ref{ThmB}(2) holds, since $\mathcal C_{\bm \lambda,\infty} \subset \mathscr{X}_p^{\infty}(\mathcal P_{m_0})$ if and only if $\lambda_1+\dots+\lambda_{m_0} = p$.
    \end{proposition}
    \begin{proof}
        Let $p_k := \lambda_1 + \dots + \lambda_k \in \overline{\mathbb Z}_{\geq 0}$, then each $\mathscr{X}_{\geq p_k}^{\infty}(\mathcal P_k) := \bigsqcup_{q\geq p_k}\mathscr{X}_{q}^{\infty}(\mathcal P_k)\sqcup (V(\mathcal P_k))_{\infty}$ is closed in $\mathcal{M}_{\infty}$. We also have the following observation:
        \begin{displaymath}
            \mathcal{C}_{\bm \lambda,\infty} = \bigcap_{k} \mathscr X_{p_k}^{\infty}(\mathcal P_k),
        \end{displaymath}
        where we write $\mathscr{X}_{\infty}^{\infty}(\mathcal P_k) = (V(\mathcal P)_k)_{\infty}$.
        
        On the one hand, suppose $\mathcal C_{\bm \lambda',\infty} \subset \overline{\mathcal C_{\bm \lambda,\infty}}$. Then we have $\mathcal C_{\bm \lambda',\infty} \subset \mathscr{X}_{\geq p_k}^{\infty}(\mathcal P_{k})$ for all $k$. Equivalently, $\sum_{i=1}^k \lambda_i \geq p_k = \sum_{i = 1}^k \lambda_i$ for all $k$. On the other hand, suppose $\bm \lambda \preceq \bm \lambda'$. Applying the same reduction as in \cite[Subsection 2.4 and Section 4]{Roi}, we may assume $n=2$, and only consider the cases: (1) $(\lambda_1,\lambda_2) = (\lambda_1'-1,\lambda_2')$; (2) $(\lambda_1,\lambda_2) =(\lambda_1'-1,\lambda_2'+1)$.  
        
        For the first case, for $s\in\mathbb C$, consider
        \begin{displaymath}
            w_1 = \begin{pmatrix}
                (st^{\lambda_1}+t^{\lambda_1+1})J \\
                & t^{\lambda_2}J
            \end{pmatrix}.
        \end{displaymath}
        If $s \neq 0$, then $w_1 \in \mathcal C_{\bm \lambda,\infty}$. If $s=0$, then $w_1 \in \mathscr X_{\lambda_1+1}(\mathcal P_1) \cap \mathscr{X}_{\lambda_2}(\mathcal P_2) = \mathcal C_{\bm \lambda',\infty}$. Hence, we have $\mathcal C_{\bm \lambda',\infty} \subset \overline{\mathcal C_{\bm \lambda,\infty}}$.

        For the second case, again for $s\in \mathbb C$, consider
        \begin{displaymath}
            w_2 =
            \begin{pmatrix}
                0 & st^{\lambda_1} & t^{\lambda_1+1} & 0 \\
                -st^{\lambda_1} & 0 & 0 & -t^{\lambda_2-1}\\
                -t^{\lambda_1+1} & 0 & 0 & st^{\lambda_2}\\
                0 & t^{\lambda_2-1} & -st^{\lambda_2} & 0
            \end{pmatrix}.
        \end{displaymath}
        One can compute to get $\mathrm{Pf}(w_2) = (s^2+1) t^{\lambda_1+\lambda_2}$. If $s\neq 0,\pm \sqrt{-1}$, then $w_2 \in \mathscr{X}_{\lambda_1}(\mathcal P_1) \cap \mathscr{X}_{\lambda_1+\lambda_2}(\mathcal P_1)$. If $s = 0$, then $w_2 \in \mathscr{X}_{\lambda_1+1}(\mathcal P_1) \cap \mathscr{X}_{\lambda_1+\lambda_2}(\mathcal P_1)$. We thus conclude the proof.
    \end{proof}
    
	
	\section{Monodromy conjecture}\label{MC}
	
	\subs{\bf Canonical log resolutions of Pfaffian ideals.}  
    The canonical log resolution of $(\mathcal M,\mathcal P_{m_0})$ is given as follows. Inductively, let $\mathcal M^0 = \mathcal M$ and $\mu_i : \mathcal M^{i} \to \mathcal M^{i-1}$ be the blow-up at the strict transform of $\mathcal P_i$ under $\eta_i := \mu_1\circ \dots \circ\mu_{i-1} : \mathcal M^i \to \mathcal M^0 = \mathcal M$. Let $E_i$ be the exceptional divisor of $\eta_i$. If $m = 2m_0$, we also denote by $E_{m_0}$ the zero locus of the strict transform of $\mathcal P_{m_0}$. Following the same argument as in \cite[5.2]{LRWW17}, one deduces that $\eta := \eta_{m_0}$ is a log resolution of $(\mathcal M,\mathcal P_{m_0})$ and we have
	\begin{equation}\label{numerical}
		\mathcal P_{m_0} \cdot \mathcal O_{\mathcal M^{m_0}} = \sum_{i = 1}^{m_0} (m_0-i+1)E_i \text{ and } K_{\eta} = \sum_{i=1}^{m_0} (\frac{(m-2i+2)(m-2i+1)}{2}-1)E_i.
	\end{equation}
    Moreover, for $j = 1,\dots,m_0$, we also have
    \begin{displaymath}
        \mathcal P_{j} \cdot \mathcal O_{\mathcal M^{m_0}} = \sum_{i = 1}^{j} (j-i+1)E_i.
    \end{displaymath}
    
	\begin{proof}[Proof of Theorem \ref{ThmA}(1)(2)]
        By \eqref{Zmot}, the poles of the motivic zeta function of $\mathcal{P}_{m_0}$ are contained in
    	\begin{displaymath}
    		\{-\frac{m(m-1)}{2m_0},-\frac{(m-2)(m-3)}{2(m_0-1)},\cdots,-\frac{(m-2m_0+2)(m-2m_0+1)}{2}        \}.
    	\end{displaymath}
        However, they coincide with poles of the topological zeta function, cf. Theorem \ref{poles_of_top_zeta}. This proves (1). (2) follows immediately from  \eqref{numerical}.
	\end{proof}    

	\subs{\bf Monodromy zeta functions of Pfaffian ideals.} Let $\eta : \mathcal M^{m_0} \to \mathcal M$ be the canonical log resolution of $\mathcal P_{m_0}$ and adopt the notations in Section \ref{Preliminaries} and Section \ref{contact_loci}. In our settings, 
    $(Y,X,\mu,\mathcal I,S) = (\mathcal M^{m_0},\mathcal M,\eta,\mathcal P_{m_0},\{1,\dots,m_0\})$.
    Our computation of the monodromy zeta functions of $\mathcal P_{m_0}$ follows the framework of \cite[Section 4]{CZ25}. 

    Fix $i \in \mathcal{S}_{m_0}$ and $e\in \widehat{E}$. Let $\bm e_i \in \mathbb Z_{\geq 0}^{m_0}$ be the $i$-th standard unit vector and $l \gg 0$ be an integer. By Theorem \ref{structure_thm_of_contact_loci}, we have 
    \begin{equation}\label{Eicirc}
        [E_i^\circ \cap \phi^{-1}(e)] = \frac{[(\varphi_{l,0}^Y)^{-1}(E_i^\circ\cap \phi^{-1}(e)) \cap \mathscr{Y}_{\bm e_i}^l]}{(\mathbb L-1)\mathbb L^{lm(m-1)/2-1}},
    \end{equation}
	and 
    \begin{equation}\label{mul}
        \mu_l : (\varphi_{l,0}^Y)^{-1}(E_i^\circ\cap \phi^{-1}(e)) \cap \mathscr{Y}_{\bm e_i}^l \to \mu_l((\varphi_{l,0}^Y)^{-1}(E_i^\circ\cap \phi^{-1}(e)) \cap \mathscr{Y}_{\bm e_i}^l)
    \end{equation}
    is a Zariski locally trivial $\mathbb A^{\nu_i-1}$-fibration. Moreover, we can characterize $(\varphi_{l,0}^Y)^{-1}(E_i^\circ\cap \phi^{-1}(e)) \cap \mathscr{Y}_{\bm e_i}^l$ in the following way:

    Let $f_1,\dots,f_{M}$ be the generators of $\mathcal{P}_{m_0}$, where $M = {m \choose 2m_0}$. Then points in $\widehat{X}$ can be expressed as $(A,[\bm b]) \in  \mathcal M \times \mathbb P^{M-1}$ such that $b_iP_j(\bm a) = b_jP_i(\bm a)$ for all $i,j$. Suppose $e = (B,[\bm \theta])$ and further assume $\theta_1 = 1$. 
    Denote by $\mathcal{D}_{l,e}\subset \mathcal{M}_l \setminus (V(\mathcal{P}_{m_0}))_l$ the set of $l$-jets satisfying (1) $\gamma(0) = B$; (2) $\mathrm{ord}_t(f_1(\gamma)) \leq \mathrm{ord}_t(f_i(\gamma))$ for all $2\leq i\leq M$; (3) $\frac{f_1(\gamma)}{f_i(\gamma)}(0) = \theta_i$ for all $2\leq i\leq M$. Then we have the decomposition
    \begin{displaymath}
        (\varphi_{l,0}^Y)^{-1}(E_i^\circ\cap \phi^{-1}(e) ) \cap \mathscr{Y}_{\bm e_i}^l = \bigsqcup_{\bm \lambda \in \Omega_{n,l}^{i,m_0}} \mu_l^{-1}(\mathcal{C}_{\bm \lambda,l}\cap \mathcal{D}_{l,e}),
    \end{displaymath}
    where
    \begin{displaymath}
        \Omega_{n,l}^{i,m_0} := \{\bm \lambda \in \Lambda_{n,l}|\lambda_1 = \dots = \lambda_{i-1} = 0,\lambda_i = \dots = \lambda_{m_0} = 1\}.
    \end{displaymath}
    Together with \eqref{Eicirc} and \eqref{mul}, we have
    \begin{equation}\label{final_Ei}
        [E_i^\circ \cap \phi^{-1}(e)] = \sum_{\bm \lambda \in \Omega_{n,l}^{i,m_0}} \frac{[\mathcal{C}_{\bm \lambda,l} \cap \mathcal{D}_{l,e}] \mathbb L^{\nu_i-1}}{(\mathbb L-1)\mathbb L^{lm(m-1)/2-1}}.
    \end{equation}
    \begin{remark}
        Since we only concern the Euler characteristic of $E_i^\circ \cap \phi^{-1}(e)$, to avoid subtlety from non-integrality of $K_0(\mathrm{Var}_{\mathbb C})$, one may specialize the expression to the virtual Poincar\'e level.
    \end{remark}

    Next, we compute each $[\mathcal{C}_{\bm \lambda,l}\cap \mathcal{D}_{l,e}]$. Consider the morphism
    \begin{equation}\label{alphalambda}
        \alpha_{\bm \lambda,l}: \mathcal C_{\bm \lambda,l} \to \widehat X,\ \gamma \mapsto (\gamma(0),[\frac{f_1(\gamma)}{t^{m_0-i+1}}(0),\dots,\frac{f_M(\gamma)}{t^{m_0-i+1}}(0)]).
    \end{equation}
    If $e$ lies in the image of $\alpha_{\bm \lambda,l}$, then $\mathcal{C}_{\bm \lambda,l}\cap \mathcal{D}_{l,e} = \alpha^{-1}_{\bm \lambda,l}(e)$. Otherwise $\mathcal{C}_{\bm \lambda,l}\cap \mathcal{D}_{l,e} = \varnothing$. Now if $\mathcal{C}_{\bm \lambda,l}\cap \mathcal{D}_{l,e} \neq \varnothing$, then there exists $\gamma \in \mathcal C_{\bm \lambda,l}$ such that $h(e) = \gamma(0)$. Hence, the rank of $h(e)$ equals $i-1$ and $E_j^\circ \cap \phi^{-1}(e) = \varnothing$ for all $j \neq i$.
    
    Note that $\alpha_{\bm \lambda,l}(\mathcal C_{\bm \lambda,l})$ can be endowed with a natural transitive $G_l = (\mathrm{GL}_m(\mathbb C))_l$-action, also induced by the congruent transformation. This makes $\alpha_{\bm \lambda,l}$ an equivariant morphism of $G_l$-varieties and it follows that
    \begin{displaymath}
        [\mathcal{C}_{\bm \lambda,l}\cap \mathcal{D}_{l,e}] = \frac{[\mathcal C_{\bm \lambda,l}]}{[\alpha_{\bm \lambda,l}(\mathcal C_{\bm \lambda,l})]}.
    \end{displaymath}
    We have to compute $[\alpha_{\bm \lambda,l}(\mathcal C_{\bm \lambda,l})]$.

    \begin{proposition}
        Let $\bm \lambda \in \Omega_{n,l}^{i,m_0}$ and $q$ be the largest integer such that $\lambda_q = 1$. If $q= m_0$, then
        \begin{displaymath}
            [\alpha_{\bm \lambda,l}(\mathcal C_{\bm \lambda,l})] = \frac{[\mathrm{GL}_m]}{[\mathrm{Sp}_{2(i-1)}] \times \mathbb L^{2(i-1)(m-2(i-1))} \times [\mathrm{GL}_{2(q-i+1)}] \times \mathbb L^{2(q-i+1)(m-2q)} \times [\mathrm{GL}_{m-2q}]}.
        \end{displaymath}
        If $q > m_0$, then
        \begin{displaymath}
            [\alpha_{\bm \lambda,l}(\mathcal C_{\bm \lambda,l})] = \frac{[\mathrm{GL}_m]}{[\mathrm{Sp}_{2(i-1)}] \times \mathbb L^{2(i-1)(m-2(i-1))} \times [\mathrm{CSP}_{2(q-i+1)}] \times \mathbb L^{2(q-i+1)(m-2q)} \times [\mathrm{GL}_{m-2q}]},
        \end{displaymath}
        where $\mathrm{CSP}_{2(q-i+1)}$ is the conformal symplectic group.
    \end{proposition}
    \begin{proof}
        Note that the $G_l$-action on $\alpha_{\bm \lambda,l}(\mathcal{C}_{\bm \lambda,l})$ factors through $G_0 = \mathrm{GL}_m(\mathbb C)$. It can be expressed as
        \begin{displaymath}
            P \cdot \alpha_{\bm \lambda,l}(\delta_{\bm \lambda,l}) = (P \begin{pmatrix}
                S_{i-1} & O\\
                O & O
            \end{pmatrix}P^T, [\frac{f_1(P\delta_{\bm \lambda}P^T)}{t^{m_0-i+1}}(0),\dots,\frac{f_M(P\delta_{\bm \lambda}P^T)}{t^{m_0-i+1}}(0)]),
        \end{displaymath}
        where $S_{i-1}$ is the block diagonal matrix given by $i-1$ copies of $J$.
        
        Denote by $Z_{\bm \lambda} \subset \mathrm{GL}_m(\mathbb C)$ the stabilizer of $\alpha_{\bm \lambda,l}(\delta_{\bm \lambda,l})$. For $P \in \mathrm{GL}_m(\mathbb C)$, $P \in Z_{\bm \lambda}$ if and only if 
        
        \noindent(1) We have $ P \begin{pmatrix}
                S_{i-1} & O\\
                O & O
            \end{pmatrix}P^T = \begin{pmatrix}
                S_{i-1} & O\\
                O & O
            \end{pmatrix}$.

        \noindent(2) There exists $a\in \mathbb C^*$, such that $$\frac{\mathrm{Pf}(G_{1,\dots,2i-2,2k_{i}-1,2k_i,\dots,2k_{m_0}-1,2k_{m_0}})}{t^{m_0-i+1}}(0) = a$$ 
        for all $i \leq k_{i} < \dots < k_{m_0} \leq q$, where $G_{1,\dots,2k_{m_0}-1,2k_{m_0}}$ is the principal submatrix of $G:= P\delta_{\bm \lambda ,l} P^T$ induced by the rows and columns of indices $1,\dots,2k_{m_0}-1,2k_{m_0}$.

        \noindent(3) The Pfaffians of other $2m_0 \times 2m_0$ submatrices of $G$ are zero.

        Suppose $P = \begin{pmatrix}
            P_1 & P_2 \\
            P_3 & P_4
        \end{pmatrix}$ and $\delta_{\bm \lambda,l} = \begin{pmatrix}
                S_{i-1} & O\\
                O & O
            \end{pmatrix} + t\begin{pmatrix}
                O & O\\
                O & M
            \end{pmatrix}$. Then (1) is equivalent to $P_1 \in \mathrm{Sp}_{2(i-1)}$, $P_3 = 0$, and $P_4$ is invertible. 
            
            Let $G' = P_4MP_4^T(0) \in M_{m-2(i-1)}^{\mathrm{skew}}(\mathbb C)$. Then (2) and (3) can be respectively rephrased as 
        
        \noindent(2$'$) There exists $a\in \mathbb C^*$, such that $$\mathrm{Pf}(G'_{2k_{i}-1,2k_i,\dots,2k_{m_0}-1,2k_{m_0}})(0) = a$$ 
        for all $i \leq k_{i} < \dots < k_{m_0} \leq q-i+1$.

        \noindent(3$'$) The Pfaffians of other $2(m_0-i+1) \times 2(m_0-i+1)$ principal submatrices of $G'$ are zero.

        By Lemma \ref{technique_computation} after this proof, one deduces from (2$'$)(3$'$) that $G'$ is of the form $\begin{pmatrix}
            G'' & O \\
            O & O
        \end{pmatrix}$, where $G'' \in {M}^{\mathrm{skew}}_{2(q-i+1)}(\mathbb C)$ is invertible. Moreover, we have $G'' \in \mathbb C^* \cdot S_{q-i+1}$ if $q > m_0$. Write $P_4 = \begin{pmatrix}
            P_5 & P_6 \\
            P_7 & P_8
        \end{pmatrix}$,
        then $P_7 = 0$, $P_5 \in \mathrm{GL}_{2(q-i+1)}(\mathbb C)$, and $P_8 \in \mathrm{GL}_{m-2q}(\mathbb C)$. Moreover, $P_5 \in \mathrm{CSP}_{2(q-i+1)}$ if $q > m_0$.

        Therefore, we conclude that:

        \noindent(i) If $q = m_0$, then
        \begin{displaymath}
            Z_{\bm \lambda} = \mathrm{Sp}_{2(i-1)} \times \mathbb A^{2(i-1)(m-2(i-1))} \times \mathrm{GL}_{2(q-i+1)} \times \mathbb A^{2(q-i+1)(m-2q)} \times \mathrm{GL}_{m-2q}.
        \end{displaymath}

        \noindent(ii) If $q > m_0$, then
         \begin{displaymath}
            Z_{\bm \lambda} = \mathrm{Sp}_{2(i-1)} \times \mathbb A^{2(i-1)(m-2(i-1))} \times \mathrm{CSP}_{2(q-i+1)} \times \mathbb A^{2(q-i+1)(m-2q)} \times \mathrm{GL}_{m-2q}.
        \end{displaymath}
        We finish the proof.
    \end{proof}

    \begin{lemma}\label{technique_computation}
        Let $H = (h_{ij}) \in M_{q}^{\mathrm{skew}}(\mathbb C)$, $1\leq r \leq k \leq \lfloor q/2\rfloor$, and $a \in \mathbb C^*$. Suppose 
        $$\mathrm{Pf}(H_{2k_1-1,2k_1,\dots,2k_r-1,2k_r}) = a$$
        for all $1\leq k_1 < \dots < k_r \leq k$ and Pfaffians of other $2r\times 2r$ principal submatrices are zero. Then
        \begin{displaymath}
            H = \begin{pmatrix}
                H' & O \\
                O & O
            \end{pmatrix}
        \end{displaymath}
        for some $H' \in M_{2k}^{\mathrm{skew}}(\mathbb C)$. Moreover, if $k > r$, then $H' \in \mathbb C^* \cdot S_{k}$, where $S_{k}$ is the block diagonal matrix given by $k$ copies of $J$. 
    \end{lemma}
    \begin{proof}
        We start with the first assertion. Without loss of generality, we may assume $q = 2k+1$. For each $1\leq i\leq 2k$, consider the expansion of the Pfaffian of the following matrix with respect to its last row.
        \begin{equation}\label{Hi1}
            H_{i} := \begin{pmatrix}
                0 & \cdots & -h_{p-2,1} & -h_{p-1,1} & -h_{p1}+h_{i1} & -h_{p1}\\
                \vdots & \ddots & \vdots & \vdots & \vdots & \vdots \\
                h_{p-2,1} & \cdots & 0 & -h_{p-1,p-2} & -h_{p,p-2}+h_{i,p-2} & -h_{p,p-2}\\
                h_{p-1,1} & \cdots & h_{p-1,p-2} & 0 & -h_{p,p-1}+h_{i,p-1} & -h_{p,p-1}\\
                h_{p1}-h_{i1} & \cdots & h_{p,p-2}-h_{i,p-2} & h_{p,p-1}-h_{i,p-1} & 0 & -h_{pi} \\
                h_{p1} & \cdots & h_{p,p-2 }& h_{p,p-1} & h_{pi} & 0\\
            \end{pmatrix}.
        \end{equation}
        Then one can check that $H_i$ is not invertible and 
        \begin{align*}
                        0 & = \mathrm{Pf}(H_i) \\
                        & = -h_{pi} \cdot \mathrm{Pf}(H_{1,\dots,p-1}) + \sum_{j = 1}^{p-1} (-1)^{j} h_{pj} \cdot \mathrm{Pf}(H_{1,\dots,\hat j, \dots,p})\\
            & = -h_{pi} \cdot \mathrm{Pf}(H_{1,\dots,p-1}),
        \end{align*}
        where $\hat j$ means removing $j$ from the sequence. Hence, $h_{pi} = 0$ for all $i$.

        Next, we prove the second assertion. We may assume $p = 2k = 2(r+1)$ and first show $H$ is a block diagonal matrix given by scalar multiples of $J$. It suffices to show $a_{1i} = a_{2i} = 0$ for all $i \geq 3$ and the claim follows by induction. Similar to \eqref{Hi1}, consider 
        \begin{equation}\label{Hi2}
            \overline {H}_{i} = 
            \begin{pmatrix}
                0 & h_{1i} & h_{13}-h_{i3} & \cdots & h_{1p}-h_{ip} \\
            -h_{1i} & 0 & h_{23} & \cdots & h_{2p} \\
            -h_{13}+h_{i3} & -h_{23} & 0 & \cdots & h_{3p} \\
            \vdots & \vdots & \vdots & \ddots & \vdots\\
            -h_{1p}+h_{ip} & -h_{2p} & -h_{3p} & \cdots & 0
            \end{pmatrix}.
        \end{equation}
        Then again,
        \begin{align*}
            0 & = \mathrm{Pf}(\overline{H}_i) \\
            & = h_{pi} \cdot \mathrm{Pf}(H_{3,\dots,p}) + \sum_{j = 3}^{p} (-1)^{j} (h_{1j}-h_{ij}) \cdot \mathrm{Pf}(H_{2,\dots,\hat j, \dots,p})\\
            & = h_{pi} \cdot \mathrm{Pf}(H_{3,\dots,p}).
        \end{align*}
        Similarly, $a_{2i} = 0$ for all $i \geq 3$.

        Finally, we show that the scalar factors in the diagonal blocks are equal. Let $\lambda_1,\dots,\lambda_k$ be these scalars. Then for all $I \subset \{1,\dots,k\}$ with $\vert I \vert = r$, we have $\prod_{i\in I} \lambda_i = a$. If $r = 1$, the assertion is trivial. Otherwise, for all $i\neq j$, pick $I \subset \{1,\dots,k\}\setminus\{i,j\}$, then $\lambda_i\prod_{l\in I} \lambda_l = \lambda_j \prod_{l\in I} \lambda_l$. Hence, $\lambda_i = \lambda_j$ and we finish the proof.
    \end{proof}

    The following theorem implies Theorem \ref{ThmA}(3).
    \begin{theorem}
        Let $\bm \lambda \in \Omega_{n,l}^{i,m_0}$ and $q,s$ be the largest integers such that $\lambda_q = 1$ and $\lambda_s \leq l$, respectively. Suppose $s = q > m_0 \geq i$ or $s = q = m_0 = i$, then for $e \in \alpha_{\bm \lambda,l}(\mathcal C_{\bm \lambda,l})$, we have $Z_{\mathcal{P}_{m_0},e}^{\mathrm{mon}}(t) = (1-t^{m_0+1-i})^{k}$ for some $k > 0$.
    \end{theorem}
    \begin{proof}
        By the discussion under \eqref{alphalambda}, we have $E_j^\circ \cap \phi^{-1}(e) = \varnothing$ for all $j\neq i$, so it suffices to compute $[E_i^\circ \cap \phi^{-1}(e)]$. Let $\bm \lambda' \in \Omega_{n,l}^{i,m_0}$ be another tuple of integers and define $q',s'$ for $\bm \lambda'$ accordingly. Recall that
        \begin{align*}
            [\mathrm{Sp}_{2k}] & = \mathbb L^{k^2}(\mathbb L^{2k}-1)(\mathbb L^{2k-2}-1)\cdots (\mathbb L^2-1),\\
            [\mathrm{CSP}_{2k}] & = \mathbb L^{k^2}(\mathbb L^{2k}-1)(\mathbb L^{2k-2}-1)\cdots (\mathbb L^2-1)(\mathbb L-1),\\
            [\mathrm{GL}_k] & = \mathbb L^{\frac{k(k-1)}{2}}(\mathbb L^k-1) \cdots (\mathbb L-1).
        \end{align*}
        Hence, the orders of $\mathbb L-1$ in $[\mathcal C_{\bm \lambda',l}]$ and $[\alpha_{\bm \lambda',l}(\mathcal C_{\bm \lambda',l})]$ are $s'$ and $\kappa_{q',m_0}$ respectively, where
        \begin{displaymath}
            \kappa_{q',m_0} = \begin{cases}
                i-1, & q' = m_0,\\
                q'-1, & q' > m_0 \geq 1.
            \end{cases}
        \end{displaymath}
        Therefore, $\chi(\frac{[\mathcal{C}_{\bm \lambda',l} \cap \mathcal{D}_{l,e}]}{\mathbb L-1}) \neq 0$ if and only if (1) $\mathcal{C}_{\bm \lambda',l} \cap \mathcal{D}_{l,e}\neq \varnothing$, and (2) $s' = q' > m_0\geq i$ or $s' = q' = m_0 = i$. Let $\Omega'$ be the set of $\bm \lambda'$ satisfying (2). Next, we show $\chi(\frac{[\mathcal{C}_{\bm \lambda',l}]}{[\alpha_{\bm \lambda',l}(\mathcal C_{\bm \lambda',l})](\mathbb L-1)}) = 1$ for all $\bm \lambda' \in \Omega'$.


        

        \noindent\textit{Case 1: $s' = q' > m_0\geq i$.}
        \begin{align*}
            & \chi(\frac{[\mathcal{C}_{\bm \lambda',l}]}{[\alpha_{\bm \lambda',l}(\mathcal C_{\bm \lambda',l})](\mathbb L-1)})\\
            & = \frac{m!}{(m-2q')! \cdot (2(q'-i+1))!! \cdot (2(i-1))!!} \cdot \frac{(2(i-1))!! \cdot (2(q'-i+1))!! \cdot (m-2q')!}{m!}\\
            & = 1.
        \end{align*}

        \noindent\textit{Case 2: $s' = q' = m_0 = i$.}
        \begin{align*}
            & \chi(\frac{[\mathcal{C}_{\bm \lambda',l}]}{[\alpha_{\bm \lambda',l}(\mathcal C_{\bm \lambda',l})](\mathbb L-1)})\\
            & = \frac{m!}{(m-2q')! \cdot (2(q'-i+1))!! \cdot (2(i-1))!!} \cdot \frac{(2(i-1))!! \cdot (2(q'-i+1))! \cdot (m-2q')!}{m!}\\
            & = 1.
        \end{align*}
        Therefore, for all $\bm \lambda' \in \Omega'$, we have 
        \begin{displaymath}
            \chi(\frac{[\mathcal{C}_{\bm \lambda',l} \cap \mathcal{D}_{l,e}]}{\mathbb L-1}) = \begin{cases}
                0, & \mathcal{C}_{\bm \lambda',l} \cap \mathcal{D}_{l,e} = \varnothing,\\
                1, & \mathcal{C}_{\bm \lambda',l} \cap \mathcal{D}_{l,e} \neq \varnothing.
            \end{cases}
        \end{displaymath}
        Now by \eqref{final_Ei} and the discussion above, we have $$\chi(E_i^\circ \cap \phi^{-1}(e)) = \sum_{\bm \lambda'\in \Omega'} \chi(\frac{[\mathcal{C}_{\bm \lambda',l} \cap \mathcal{D}_{l,e}]}{\mathbb L-1})\geq \chi(\frac{[\mathcal{C}_{\bm \lambda,l} \cap \mathcal{D}_{l,e}]}{\mathbb L-1}) = 1.$$ We thus conclude the proof.
    \end{proof}

    \par $\space$

    \noindent \bf{Acknowledgement.} \rm Y. Chen, Q. Shi, and H. Zuo were supported by BJNSF Grant 1252009. H. Zuo was supported by NSFC Grant 12271280.


\begin{thebibliography}{BBBBB25}

    \bibitem[AKMW02]{AKW}
	D. Abramovich, K. Karu, K. Matsuki, and J. W\l odarczyk.
	\newblock Torification and factorization of birational maps.
	\newblock {\em J. Amer. Math. Soc.}, 15(3): 531--572, 2002.
    
    \bibitem[Bla24]{Bla24}
    G. Blanco.
    \newblock{Topological roots of the Bernstein-Sato polynomial of plane curves}.
    \newblock{\em \url{https://arxiv.org/abs/2406.09034}}, 2024.
    
        \bibitem[BBPFP24]{BBPP24}
		N. Budur, J. de la Bodega, E. de Lorenzo Poza, J. Fern\'andez de Bobadilla, and T. Pe\l ka.
		\newblock{On the embedded Nash problem}.
		\newblock{\em Forum Math. Pi}, 12, Paper No. e15, 28 pp., 2024.
            
        \bibitem[BFLN22]{Cohomology_of_Contact_Loci}
		N. Budur, J. Fern\'andez de~Bobadilla, Q.~T. L\^e, and H.~D.
		Nguyen.
		\newblock Cohomology of contact loci.
		\newblock {\em J. Differential Geom.}, 120(3): 389--409, 2022.
    

    
		

        \bibitem[BMS06]{BMS06}
        N.~Budur, M.~Musta\c{t}\v{a}, and M.~Saito.
        \newblock Bernstein-{S}ato polynomials of arbitrary varieties.
        \newblock {\em Compos. Math.}, 142(3): 779--797, 2006.

       \bibitem[BT20]{BT20}
        N. Budur and T. Q. Tue.
        \newblock{On contact loci of hyperplane arrangements.}
        \newblock{\em Adv. in Appl. Math.}, 132, 102271, 2022.

        \bibitem[BSZ25]{BSZ25}
        N. Budur, Q. Shi, and H. Zuo.
        \newblock {Polar loci of multivariable archimedean zeta functions.}
        \newblock {\em \url{https://arxiv.org/abs/2504.10051}}, 2025.
    
        \bibitem[BVWZ21]{BVWZ21}
        N. Budur, R. van~der Veer, L. Wu, and P. Zhou.
        \newblock {Zero loci of {B}ernstein-{S}ato ideals.}
        \newblock {\em Invent. Math.}, 225: 45-72, 2021.

        \bibitem[CIL17]{CIL17}
        W. Castryck, D. Ibadula, and A. Lemahieu.
        \newblock{The holomorphy conjecture for nondegenerate surface singularities.}
        \newblock{\em Nagoya Math. J.}, 227: 160--188, 2017.

    
        \bibitem[CLNS18]{CLNS}
        A. Chambert-Loir, J. Nicaise, and J. Sebag. 
        \newblock{Motivic integration}.
        \newblock{volume 325 of Progress in Mathematics.} Birkhäuser/Springer, New York, 2018.
    
        \bibitem[CZ25]{CZ25}
        Y. Chen and H. Zuo.
        \newblock{On the monodromy conjecture for determinantal varieties}.
        \newblock{\em \url{https://arxiv.org/abs/2510.11425}}, 2025.

        \bibitem[FD16]{FD16}
        T. de~Fernex and R. Docampo.
        \newblock{Terminal valuations and the Nash problem}.
        \newblock{Invent. Math.}, 203(1): 303--331, 2016.

        \bibitem[Den91]{Den}
        J. Denef.
        \newblock{Report on Igusa's local zeta function.}
        \newblock{\em Ast\'erisque}, 201-203: 359--386, 1991.
        
		\bibitem[DL92]{DL92}
				J.~Denef and F.~Loeser.
				\newblock Caract\'{e}ristiques d'{E}uler-{P}oincar\'{e}, fonctions z\^{e}ta locales et modifications analytiques.
				\newblock {\em J. Amer. Math. Soc.}, 5(4):705--720, 1992.


        
		\bibitem[DL98]{DL98}
		J. Denef and F. Loeser.
		\newblock Motivic {I}gusa zeta functions.
		\newblock {\em J. Algebraic Geom.} 7: 505-537, 1998.
		
		\bibitem[Doc13]{Roi}
		R. Docampo. 
		\newblock{Arcs on determinantal varieties.}
		\newblock{\em Trans. Amer. Math. Soc.}, 365: 2241-2269, 2013.


        
		\bibitem[ELM04]{ELM04}
		L. Ein, R. Lazarsfeld, and M. Musta\c t\u a.
		\newblock{Contact loci in arc spaces}.
		\newblock{\em Compos. Math.}, {140}(5): 1229--1244, 2004.		

        \bibitem[ELT22]{E+}
        A. Esterov, A. Lemahieu, and K. Takeuchi.
        \newblock{On the monodromy conjecture for non-degenerate hypersurfaces},
        \newblock{\em J. Eur. Math. Soc.}, 24(11): 3873--3949, 2022.
        

    
        
		\bibitem[Igu00]{Igu} J.-I. Igusa. {An introduction to the theory of local zeta functions.}
		AMS/IP Stud. Adv. Math., 14,  xii+232 pp., 2000.

        \bibitem[Ish04]{Ish04}
        S. Ishii.
        \newblock{The arc space of a toric variety.}
        \newblock{\em J. Algebra}, 278: 666-683, 2004.

        \bibitem[IK03]{IK03}
        S. Ishii and J. Koll\'ar.
        \newblock{The Nash problem on arc families of singularities.}
        \newblock{\em Duke Math. J.}, 120(3): 601--620, 2003.
            
		\bibitem[LRWW17]{LRWW17}
		A. L\H{o}rincz, C. Raicu, U. Walther, and J. Weyman.
		\newblock{Bernstein-{S}ato polynomials for maximal minors and
			sub-maximal {P}faffians}.
		\newblock{\em Adv. Math.}, 307: 224--252, 2017.

        \bibitem[Loe88]{Loe88}
        F. Loeser.
        \newblock{Fonctions d’Igusa p-adiques et polyn\^{o}mes de Bernstein}
        \newblock{\em Amer. J. Math.}, 110(1): 1--21, 1988. 
		
%
%
%
%
%
%
%

%
%
%

%
%
%
%
%
%
%
%
%
%
%
%
%
%
%
%
%
%
%
            \bibitem[LVP11]{VPL11}
            A. Lemahieu and L. Van Proeyen.
            \newblock{The holomorphy conjecture for ideals in dimension two}.
            \newblock{\em Proc. Amer. Math. Soc.}, 139(11): 3845--3852, 2011.

            \bibitem[RV01]{RV01}
            B. Rodrigues and W. Veys.
            \newblock{Holomorphy of Igusa's and topological zeta functions for homogeneous polynomials.}
            \newblock{\em Pacific J. Math.}, 201(2): 429--440, 2001.

            \bibitem[SZ24]{SZ24}
            Q. Shi and H. Zuo.
            \newblock{On the Tensor Property of Bernstein-Sato Polynomial}.
            \newblock{\em \url{https://arxiv.org/abs/2406.04121}}, 2024.
            
            \bibitem[VPV10]{PV10}
    		  L. Van Proeyen and W. Veys. 
    		  \newblock{The monodromy conjecture for zeta functions associated to ideals in dimension two}.
    		\newblock{\em Ann. Inst. Fourier (Grenoble)}, 60(4): 1347-1362, 2010.
        
			\bibitem[Ver83]{Ver83}
            J.-L. Verdier. 
            \newblock{Spécialisation de faisceaux et monodromie modérée.}
            \newblock{\em Astérisque}, 101-102: 332-364, 1983.

            \bibitem[Vey90]{Vey90}
            W. Veys.
            \newblock{On the poles of Igusa’s local zeta function for curves.}
            \newblock{\em J. Lond. Math. Soc.}, 41(1): 27--32, 1990.

            \bibitem[Vey93]{Vey93}
            W. Veys. Holomorphy of local zeta functions for curves.
            \newblock{\em Math. Ann.}, 295: 635--641, 1993.   

            
            \bibitem[Vey06]{Vey06}
            W.~Veys.
            \newblock Vanishing of principal value integrals on surfaces.
            \newblock {\em J. Reine Angew. Math.}, 598: 139--158, 2006.
            
            \bibitem[Vey25]{Vey25}
			W. Veys.
			\newblock Introduction to the monodromy conjecture.
            \newblock{\em Handbook of geometry and topology of singularities VII}, Springer, Cham, 721–765, 2025.
%
		
		
	\end{thebibliography}
\end{document}